\documentclass[12pt]{article}

\usepackage{amsmath}
\usepackage{amssymb}
\usepackage{amsthm}
\usepackage[pdfpagemode=None]{hyperref}
\usepackage{graphicx}
\usepackage{pstricks}

\newtheorem{thm}{Theorem}[section]
\newtheorem{co}[thm]{Corollary}
\newtheorem{lem}[thm]{Lemma}

\newtheorem{assumption}[thm]{Assumption}

\newtheorem{pr}[thm]{Proposition}

\newtheorem{definition}[thm]{Definition}
\newenvironment{de}{\begin{definition}\rm}{\end{definition}}
\newtheorem{example}[thm]{Example}
\newenvironment{exmp}{\begin{example}\rm}{\end{example}}
\newtheorem{remark}[thm]{Remark}
\newenvironment{rem}{\begin{remark}\rm}{\end{remark}}
\newtheorem{tab}{Table}

\newcommand{\diag}{{\rm diag}\,}

\def\eps{\varepsilon}

\title{Analyticity of Entropy Rate of Hidden Markov Chains}

\author{Guangyue Han, Brian Marcus\\
  {\normalsize Department of Mathematics}\vspace{-1mm} \\
  {\normalsize University of British Columbia}\vspace{-1mm} \\
  {\normalsize Vancouver, B.C.,  V6T 1Z2}\\
  {\normalsize {\em e-mail:\/} {ghan, marcus}@math.ubc.ca}}

\date{{\normalsize \today}}

\begin{document}\maketitle\thispagestyle{empty}


\begin{abstract}
We prove that under mild positivity assumptions the entropy rate of
a hidden Markov chain varies analytically as a function of the
underlying Markov chain parameters. A general principle to determine
the domain of analyticity is stated. An example is given to estimate
the radius of convergence for the entropy rate. We then show that
the positivity assumptions can be relaxed, and examples are given
for the relaxed conditions. We study a special class of hidden
Markov chains in more detail: binary hidden Markov chains with an
unambiguous symbol, and we give necessary and sufficient conditions
for analyticity of the entropy rate for this case. Finally, we show
that under the positivity assumptions the hidden Markov chain {\em
itself} varies analytically, in a strong sense, as a function of the
underlying Markov chain parameters.
\end{abstract}

\section{Introduction}
For $m, n \in \mathbb{Z}$ with $m \leq n$,  we denote a sequence of
symbols $y_m, y_{m+1}, \ldots, y_n$ by $y_m^n$. Consider a
stationary stochastic process $Y$ with a finite set of states $\{1,
2, \cdots, B\}$ and distribution $p(y_m^n)$.  Denote the conditional
distributions by $p(y_{n+1}|y_m^n)$. The entropy rate of $Y$ is
defined as
$$
H(Y)=\lim_{n \to \infty} -E_p(\log(p(y_0|y_{-n}^{-1}))),
$$
where $E_p$ denotes expectation with respect to the distribution
$p$.

Let $Y$ be a stationary first order Markov chain with
$$
\qquad \Delta(i, j)=p(y_1=j | y_0=i).
$$
It is well known that
$$
H(Y)=-\sum_{i, j} p(y_0=i) \Delta(i, j) \log \Delta(i, j).
$$

A \emph{hidden Markov chain} $Z$ (or function of a Markov chain) is
a process of the form $Z=\Phi(Y)$, where $\Phi$ is a function
defined on $\{1, 2, \cdots, B\}$ with values $\{1, 2, \cdots, A\}$.
Often a hidden Markov chain is defined as a Markov chain observed in
noise. It is well known that the two definitions are equivalent (the
equivalence is typified by Example~\ref{bsc}).

For a hidden Markov chain, $H(Z)$ turns out   (see
Equation~(\ref{blackwell_form}) below) to be the integral of a
certain function defined on a simplex with respect to a measure due
to Blackwell~\cite{bl57}. However Blackwell's measure is somewhat
complicated and the integral formula appears to be difficult to
evaluate in most cases.

Recently there has been a rebirth of interest in computing the
entropy rate of a hidden Markov chain, and many approaches have been
adopted to tackle this problem.  For instance, some researchers have
used Blackwell's measure to bound the entropy rate~\cite{or03} and
others introduced a variation~\cite{eg04} on bounds due
to~\cite{bi62}.

In a new direction,~\cite{or03,ja04,zu04} have studied the variation
of the entropy rate as parameters of the underlying Markov chain
vary. These works motivated us to consider the general question of
whether the entropy rate of a hidden Markov chain is smooth, or even
analytic~\cite{sh92, ta02}, as a function of the underlying
parameters. Indeed, this is true under mild positivity assumptions:

\begin{thm}  \label{main}
Suppose that the entries of $\Delta$ are analytically parameterized
by a real variable vector $\vec{\varepsilon}$. If at
$\vec{\varepsilon}=\vec{\varepsilon}_0$,
\begin{enumerate}
\item For all $a$, there is at least one $j$ with $\Phi(j)=a$ such that the
$j$-th column of $\Delta$ is strictly positive -- and --
\item Every column of $\Delta$ is either all zero or strictly
positive,
\end{enumerate}
then $H(Z)$ is a real analytic function of $\vec{\varepsilon}$ at
$\vec{\varepsilon}_0$.
\end{thm}

Note that this theorem holds if all the entries of $\Delta$ are
positive.  The more general form of our hypotheses is very important
(see Example~\ref{bsc}).

Real analyticity at a point is important because it means that the
function can be expressed as a convergent power series in a
neighborhood of the point. The power series can be used to
approximate or estimate the function. For convenience of the reader,
we recall some basic concepts of analyticity in
Section~\ref{XXX-XXX}.

Several authors have observed that the entropy rate of a hidden
Markov chain can be viewed as the top Lyapunov exponent of a random
matrix product~\cite{ho03,ja04,gh95u}.  Results in ~\cite{ar94a,
pe90, pe92, ru79} show that under certain conditions the top
Lyapunov exponent of a random matrix product varies analytically as
either the underlying Markov process varies analytically or as the
matrix entries vary analytically, but not both. However, when
regarding the entropy rate as a Lyapunov exponent of a random matrix
product, the matrix entries depend on the underlying Markov process.
So, the results from Lyapunov theory do not appear to apply
directly. Nevertheless, much of the main idea of our proof of
Theorem~\ref{main} is essentially contained in Peres~\cite{pe92}. In
contrast to Peres' proof, we do not use the language of Lyapunov
exponents and we use only basic complex analysis and no functional
analysis. Also the hypotheses in~\cite{pe92} do not carry over to
our setting. To the best of our knowledge the statement and proof of
Theorem~\ref{main} has not appeared in the literature. For
analyticity of certain other statistical quantities, see also
related work in the area of statistical physics in~\cite{do73, ca81,
lo98, ch03}.

After discussing background in Sections~\ref{xxx-III}
and~\ref{XXX-XXX}, we prove Theorem~\ref{main} in
Section~\ref{xxx-IV}.  As an example, we show that the entropy rate
of a hidden Markov chain obtained by observing a binary Markov
chains in binary symmetric noise, with noise parameter $\eps$, is
analytic at any $\eps = \eps_0 \ge 0$, provided that the Markov
transition probabilities are all positive.

In Section~\ref{domain}, we infer from the proof of
Theorem~\ref{main} a general principle to determine a domain of
analyticity for the entropy rate. We apply this to the case of
hidden Markov chains obtained from binary Markov chains in binary
symmetric noise to find a lower bound on the radius of convergence
of a power series in $\eps$ at $\eps_0 = 0$.  Given the recent
results of~\cite{zu05}, which compute the derivatives of all orders
at $\eps_0=0$, this gives an explicit power series for entropy rate
near $\eps_0=0$.

In Section~\ref{relaxed}, we
show how to relax the conditions of Theorem~\ref{main}
and apply this to give more examples where the entropy rate
is analytic.

The entropy rate can fail to be analytic.
In Section~\ref{boundary} we give examples and then give a complete
set of necessary and sufficient conditions for analyticity in the
special case of binary hidden Markov chains with
an unambiguous symbol, i.e., a symbol which can be produced by only
one symbol of the Markov chain.

Finally in Section~\ref{xxx-II}, we resort to more advanced
techniques to prove a stronger version, Theorem~\ref{main-1}, of
Theorem~\ref{main}. This result gives a sense in which the hidden
Markov chain {\em itself} varies analytically with
$\vec{\varepsilon}$. The proof of this result requires some measure
theory and functional analysis, along with ideas from equilibrium
states~\cite{ru78}, which are reviewed in
Appendix~\ref{equilibrium}. Our first proof of Theorem~\ref{main}
was derived as a consequence of Theorem~\ref{main-1}. It also
follows from Theorem~\ref{main-1} that, in principle, many
statistical properties in addition to entropy rate vary
analytically.

Most results of this paper were first announced in~\cite{gm05}.

\section{Iteration on the Simplex} \label{xxx-III}

Let $W$ be the simplex, comprising the vectors
$$
\{w=(w_1, w_2, \cdots, w_B) \in \mathbb{R}^B:w_i \geq 0, \sum_i
w_i=1\},
$$
and let $W_a$ be all $w \in W$ with $w_i=0$ for $\Phi(i) \neq a$.
Let $W^{\mathbb{C}}$ denote the complex version of $W$, i.e.,
$W^{\mathbb{C}}$ denotes the complex simplex comprising the vectors
$$
\{w=(w_1, w_2, \cdots, w_B) \in \mathbb{C}^B : \sum_i w_i=1\},
$$
and let $W_a^{\mathbb{C}}$ denote the complex version of $W_a$,
i.e., $W_a^{\mathbb{C}}$ consists of all $w \in W^{\mathbb{C}}$ with
$w_i=0$ for $\Phi(i) \neq a$. For $a \in A$, let $\Delta_a$ denote
the $B \times B$ matrix such that $\Delta_a(i, j)=\Delta(i, j)$ for
$j$ with $\Phi(j)=a$, and $\Delta_a(i, j)=0$ otherwise. For $a \in
A$, define the scalar-valued and vector-valued functions $r_a$ and
$f_a$ on $W$ by
$$
r_a(w)= w \Delta_a \mathbf{1},
$$
and
$$
f_a(w)=w \Delta_a/ r_a(w).
$$
Note that $f_a$ defines the action of the matrix $\Delta_a$ on the
simplex $W$. For any fixed $n$ and $z_{-n}^0$, define
\begin{equation} \label{x-i}
x_i=x_i(z_{-n}^i)=p(y_i=\cdot \;|z_i, z_{i-1}, \cdots, z_{-n}),
\end{equation}
(here $\cdot$ represent the states of the Markov chain $Y$,) then
from Blackwell~\cite{bl57}, $\{x_i\}$ satisfies the random dynamical
iteration
\begin{equation}
\label{iter0} x_{i+1}=f_{z_{i+1}}(x_i),
\end{equation}
starting with
\begin{equation}
\label{init0} x_{-n-1} = p(y_{-n-1}=\cdot\;).
\end{equation}

We remark that Blackwell showed that
\begin{equation}
\label{blackwell_form} H(Z)=-\int \sum_a r_a(w) \log r_a(w) dQ(w) ,
\end{equation}
where $Q$, known as {\em Blackwell's measure}, is the limiting
probability distribution, as $n \rightarrow \infty$, of $\{x_0\}$ on
$W$. However, we do not use Blackwell's measure explicitly in this
paper.

%
%

Next, we consider two metrics on a compact subset $S$ of the
interior of a subsimplex $W'$ of $W$.  Without loss of generality, we assume that $W'$
consists of all points from $W$ with the last $B-k$ coordinates
equal to $0$. The Euclidean metric $d_\textbf{E}$ on $S$ is defined
as usual, namely for $u,v \in S$,
$$
u=(u_1,u_2,\cdots,u_B), v=(v_1,v_2,\cdots,v_B) \in S,
$$
we have
$$
d_\textbf{E}(u,
v)=\sqrt{(u_1-v_1)^2+(u_2-v_2)^2+\cdots+(u_k-v_k)^2}.
$$
The Hilbert metric~\cite{se80} $d_\textbf{B}$ on $S$ is defined as
follows:
$$
d_\textbf{B}(u, v)=\max_{i \neq j \leq k} \log \left(
\frac{u_i/u_j}{v_i/v_j} \right).
$$

The following result is well known  (for instance, see
\cite{ar94a}). For completeness, we give a detailed proof in
Appendix~\ref{Metric-Equivalence1}.

\begin{pr}  \label{Metric-Equivalence}
$d_\textbf{E}$ and $d_\textbf{B}$ are equivalent (denoted by
$d_\textbf{E} \sim d_\textbf{B}$) on any compact subset $S$ of the
interior of a subsimplex $W'$ of $W$, i.e., there are positive constants $C_1 < C_2$
such that for any two points $u, v \in S$,
$$
C_1 d_\textbf{B}(u, v) < d_\textbf{E}(u, v) < C_2 d_\textbf{B}(u,
v).
$$
\end{pr}

\begin{pr} \label{eventually-contract}
Assume that at $\vec{\varepsilon}_0$, $\Delta$ satisfies conditions
$1$ and $2$ of Theorem~\ref{main}. Then for sufficiently large $n$
and all choices of $a_1, \ldots, a_n$ and $b$, the mapping $f_{a_n}
\circ f_{a_{n-1}} \circ \cdots \circ f_{a_1}$ is a contraction
mapping under the Euclidean metric on $W_b$.
\end{pr}

\begin{proof}
$\hat{W}_b=f_b(W)$ is a compact subset of the interior of some
subsimplex $W_b'$ of $W_b$; this subsimplex corresponds to column
indices $j$ such that $\Phi(j)=b$ and the $j$-th column is strictly
positive. Therefore one can define the Hilbert metric accordingly on
$\hat{W}_b$. Each $f_a$ is a contraction mapping on each $\hat{W}_b$
under the Hilbert metric~\cite{se80}; namely there exists $0 < \rho
< 1$ such that for any $a$ and $b$, and for any two points $u, v \in
\hat W_b$,
$$
d_\textbf{B}(f_a(u), f_a(v)) < \rho d_\textbf{B}(u, v).
$$
Thus, for any choices of $a_2, a_3, \cdots, a_n$, we have
$$
d_\textbf{B}(f_{a_n} \circ f_{a_{n-1}} \circ \cdots \circ f_{a_2}
(u), f_{a_n} \circ f_{a_{n-1}} \circ \cdots \circ f_{a_2} (v)) <
\rho^{n-1} d_\textbf{B}(u, v).
$$
By Proposition~\ref{Metric-Equivalence}, there exists a positive
constant $C$ such that
$$
d_\textbf{E}(f_{a_n} \circ f_{a_{n-1}} \circ \cdots \circ f_{a_2}
(u), f_{a_n} \circ f_{a_{n-1}} \circ \cdots \circ f_{a_2} (v)) < C
\rho^{n-1} d_\textbf{E}(u, v).
$$
Let $L$ be a universal Lipschitz constant for any $f_c: W_b
\rightarrow W'_c$ with respect to the Euclidean metric. Choose $n$
large enough such that $C \rho^{n-1} < 1/L$. So, for sufficiently
large $n$, any composition of the form $f_{a_n} \circ \cdots \circ
f_{a_1}$ is a Euclidean contraction on $W_b$.

\end{proof}

\begin{rem}
Using a slightly modified proof, one can show that for sufficiently
large $n$, any composition of the form $f_{a_n} \circ \cdots \circ
f_{a_1}$ is a Euclidean contraction on the whole simplex $W$.
\end{rem}

\section{Brief background on analyticity}
\label{XXX-XXX}

In this section, we briefly review the basics in complex analysis
for the purpose of this paper. For more details, we refer
to~\cite{sh92, ta02}.

A real (or complex) function of several variables is analytic at a
given point if it admits a convergent Taylor series representation
in a real (or complex) neighborhood of the given point. We say that
it is real (or complex) analytic in a neighborhood if it is real (or
complex) analytic at each point of the neighborhood.

The relationship between real and complex analytic functions is as
follows: 1)  Any real analytic function can be extended to a complex
analytic function on some complex neighborhood; 2)  Any real
function obtained by restricting a complex analytic function from a
complex neighborhood to a real neighborhood is a real analytic
function.

The main fact regarding analytic functions used in this paper is
that the uniform limit of a sequence of complex analytic functions
on a fixed complex neighborhood is complex analytic.  The analogous
statement does not hold (in fact, fails dramatically!) for real
analytic functions.



As an example of a real-valued parametrization of a matrix,
consider:
$$
\Delta(\eps)=\left [ \begin{array}{ccc}
                2\eps& \eps & 1 - 3\eps\\
                \eps & 1 - \eps - \sin(\eps) & \sin(\eps)\\
                \eps^2& \eps^3 & 1 - \eps^2 - \eps^3\\
                \end{array} \right ].
$$
Denote the states of $\Delta$ by $\{1,2,3\}$ and let
$\Phi(1)=\Phi(2) =0,~ \Phi(3) = 1$. Each entry of $\Delta$ is a real
analytic function of $\eps$ at any given point $\eps = \eps_0$. For
$\eps_0 > 0$ and sufficiently small, $\Delta$ is stochastic (i.e.,
each row sums to 1 and each entry is nonnegative) and in fact
strictly positive (i.e., each entry is positive).  According to
Theorem~\ref{main}, for such values of $\eps_0$, the entropy rate of
the hidden Markov chain defined by $\Delta(\eps)$ and $\Phi$ is real
analytic as a function of $\eps$ at $\eps_0$. .

%
%
%
%

While we typically think of analytic parametrizations as having the
``look'' of the preceding example, there is a conceptually simpler
parametrization -- namely, parameterize an $n \times n$  matrix
$\Delta$ by its entries themselves;  if $\Delta$ is required to be
stochastic, we choose the parameters to be any set of $n-1$ entries
in each row (so, the real variable vector is an $n(n-1) $-tuple).
Clearly, for analyticity it does not matter which entries are
chosen.  We call this the {\em natural parametrization}.

Suppose that $H(Z)$ is analytic with respect to this
parametrization.  Then, $H(Z)$ viewed as a function of any other
analytic parametrization of the entries of $\Delta$ is the
composition of two analytic functions and thus must be analytic. We
thus have that the following two statements are equivalent.
\begin{itemize}
\item $H(Z)$ is analytic with respect to the natural
parameterization.
\item $H(Z)$ is analytic with respect to any analytic
parameterization.
\end{itemize}
We shall use this implicitly through the paper.

\section{Proof of Theorem~\ref{main}} \label{xxx-IV}

{\em Notation:}  We rewrite $\Delta$, $Z$, $f_a(x)$,
$p(z_0|z_{-\infty}^{-1})$ with parameter vector $\vec{\varepsilon}$
as $\Delta^{\vec{\varepsilon}}$, $Z^{\vec{\varepsilon}}$,
$f_a^{\vec{\varepsilon}}(x)$ and
$p^{\vec{\varepsilon}}(z_0|z_{-\infty}^{-1})$, respectively. We use
the notation $\hat W_a$ to mean $f_a^{\eps_0}(W)$.  Let
$\Omega_{\mathbb{C}}=\Omega_{\mathbb{C}}(r)$ denote the set of
points of distance at most $r$ from $\vec{\varepsilon}_0$ in the
complex parameter space $\mathbb{C}^m$.  Let $N_b = N_b(R)$ denote
the set of all points in $W_b^{\mathbb{C}}$ of distance at most $R$
from $\hat W_b$.

We first prove that for some $r>0$, $\log
p^{\vec{\varepsilon}}(z_0|z_{-\infty}^{-1})$ can be extended to a
complex analytic function  of $\vec{\varepsilon} \in
\Omega_{\mathbb{C}}(r)$ and that $ |\log
p^{\vec{\varepsilon}}(z_0|z_{-\infty}^{-1}) - \log
p^{\vec{\varepsilon}}(\hat z_0|\hat z_{-\infty}^{-1})|$ decays
exponentially fast in $n$, when $z_{-n}^0 = \hat z_{-n}^0$,
uniformly in $\vec{\varepsilon} \in \Omega_{\mathbb{C}}(r)$.

Note that for each $a,b$, $f^{\vec{\varepsilon}}_a(w)$ is a rational
function of the entries of $\Delta^{\vec{\varepsilon}}$ and $w \in
\hat{W}_b$. So, by viewing the real vector variables
$\vec{\varepsilon}$ and $w$ as complex vector variables, we can
naturally extend $f^{\vec{\varepsilon}}_a(w)$ to a complex-valued
function of complex vector variables $\vec{\varepsilon}$  and $w$.
Since $\Delta$ satisfies conditions $1$ and $2$ at
$\vec{\varepsilon}_0$, for sufficiently small $r$ and $R$, the
denominator of $f^{\vec{\varepsilon}}_a(w)$ is nonzero for
$\vec{\varepsilon}$ in $\Omega_\mathbb{C}(r)$ and $w$ in $N_b(R)$.
Thus, $f^{\vec{\varepsilon}}_a(w)$ is a complex analytic function of
$(\vec{\varepsilon}, w)$ in the neighborhood $
\Omega_{\mathbb{C}}(r) \times N_b(R)$.

Assuming conditions $1$ and $2$, we claim that $\Delta$ has an
isolated (in modulus) maximum eigenvalue $1$ at
$\vec{\varepsilon}_0$. To see this, we apply Perron-Frobenius
theory~\cite{se80} as follows. By permuting the indices, we can
express:
$$
\Delta =
\left[
\begin{array}{cc} U & 0 \\
V & 0 \end{array} \right]
$$
where $U$ is the submatrix corresponding to indices with positive
columns.  The nonzero eigenvalues of $\Delta$ are the same as the
eigenvalues of $U$, which is a positive stochastic matrix.  Such a
matrix has isolated (in modulus) maximum eigenvalue $1$.

The stationary distribution $p^{\vec{\varepsilon}}(y=\cdot\;)$ (the
eigenvector corresponding to the maximum eigenvalue $1$) is a
rational function of the entries of $\Delta^{\vec{\varepsilon}}$,
since it is a solution of the equation $v\Delta^{\vec{\varepsilon}}
= v$. So, in the same way as for $f^{\vec{\varepsilon}}_a(w)$ we can
naturally extend $p^{\vec{\varepsilon}}(y=\cdot\;)$ to a complex
analytic function $p^{\vec{\varepsilon}}(y=\cdot\;)$ on
$\Omega_{\mathbb{C}}$.

Extending (\ref{x-i}) for each $i$, we define
\begin{equation}
\label{x-i-0}
x_i^{\vec{\varepsilon}}=x_i^{\vec{\varepsilon}}(z_{-n}^i)=p^{\vec{\varepsilon}}(y_i=\cdot
\; |z_{-n}^i),
\end{equation}
by iterating the following complexified random dynamical system
(extending (\ref{iter0}) and (\ref{init0})):
\begin{equation}
\label{iter}
x^{\vec{\varepsilon}}_{i+1}=f^{\vec{\varepsilon}}_{z_{i+1}}(x_i^{\vec{\varepsilon}}),
\end{equation}
starting with
\begin{equation}
\label{init} x^{\vec{\varepsilon}}_{-n-1} =
p^{\vec{\varepsilon}}(y_{-n-1}=\cdot \;).
\end{equation}

By Proposition~\ref{eventually-contract}, for sufficiently large
$n$, we can replace the set of mappings $\{f_a^{\eps_0}\}$ with the
set $\{f^{\eps_0}_{a_n} \circ f^{\eps_0}_{a_{n-1}} \circ \cdots
\circ f^{\eps_0}_{a_1}\}$ and then assume that each $f^{\eps_0}_a$
is a Euclidean contraction on each $W_b$ with contraction
coefficient $\rho < 1$.  Since $\hat W_b$ is compact and the
definition of $\rho$-contraction is given by strict inequality, we
can choose $r$ and $R$ sufficiently small such that
\begin{equation}
\label{contract-1}
f^{\vec{\varepsilon}}_a \mbox{ is a Euclidean } \rho- \mbox{contraction on
each } N_b(R), ~\eps \in \Omega_{\mathbb{C}}(r).
\end{equation}
Further, we claim that by choosing $r$ still smaller, if necessary,
\begin{equation}
\label{confine-orbit} x_i^{\vec{\eps}} \in \cup_b N_b(R), \mbox{ for
all } i,n \mbox{ and all choices of } z_{-n}^i,~\varepsilon \in
\Omega_{\mathbb{C}}(r).
\end{equation}
To see this, fixing $\rho$ and $R$, 
choose $r$ so small that
\begin{equation}
\label{bounded-1}
|f_{a}^{\vec{\varepsilon}}(x)-f_{a}^{\vec{\varepsilon}_0}(x)| \le
R(1-\rho),~ x \in \cup_b \hat{W}_b,~ \eps \in \Omega_{\mathbb{C}}(r)
\end{equation}
and
\begin{equation}
\label{bounded-2}
|p^{\vec{\varepsilon}}(\cdot)-p^{\vec{\varepsilon}_0}(\cdot) |\le R(1-\rho),
~\eps \in \Omega_{\mathbb{C}}(r).
\end{equation}

Now consider the difference
$$
x_{i+1}^{\vec{\varepsilon}}-x_{i+1}^{\vec{\varepsilon}_0}
$$
\begin{equation}  \label{star}
=f_{z_{i+1}}^{\vec{\varepsilon}}(x_i^{\vec{\varepsilon}})-f_{z_{i+1}}^{\vec{\varepsilon}_0}(x_i^{\vec{\varepsilon}_0})
=f_{z_{i+1}}^{\vec{\varepsilon}}(x_i^{\vec{\varepsilon}})-f_{z_{i+1}}^{\vec{\varepsilon}}(x_i^{\vec{\varepsilon}_0})+f_{z_{i+1}}^{\vec{\varepsilon}}(x_i^{\vec{\varepsilon}_0})-f_{z_{i+1}}^{\vec{\varepsilon}_0}(x_i^{\vec{\varepsilon}_0}).
\end{equation}
Then by (\ref{contract-1}) , (\ref{bounded-1}) and
(\ref{bounded-2}), and (\ref{star}),
for $i > -n-1$, we have
$$
|x_{i+1}^{\vec{\varepsilon}}-x_{i+1}^{\vec{\varepsilon}_0}| \leq
\rho |x_i^{\vec{\varepsilon}}-x_i^{\vec{\varepsilon}_0}|+R(1-\rho).
$$
So,
$$
|x_{i+1}^{\vec{\varepsilon}}-x_{i+1}^{\vec{\varepsilon}_0}| \leq R,
$$
and thus for all $i$, we have $x_{i+1}^{\vec{\varepsilon}} \in
\cup_b N_b(R)$, yielding (\ref{confine-orbit}). Each $x_i^\eps$ is
the composition of analytic functions on $\Omega_{\mathbb{C}}(r)$
and so is complex analytic on $\Omega_{\mathbb{C}}(r)$.

%

For $0 \leq n_1, n_2 \leq \infty$, we say two sequences
$\{z_{-n_1}^0\}$ and $\{\hat{z}_{-n_2}^0\}$ have a common tail if
there exists $n \geq 0$ with $n \leq n_1, n_2$ such that $
z_i=\hat{z}_i,   -n \leq i \leq 0$ (denoted by $z_{-n_1}^0
\stackrel{n}{\sim} \hat{z}_{-n_2}^0$).

Let
$$
x_i^{\vec{\varepsilon}}=x_i^{\vec{\varepsilon}}(z_{-n_1}^i)=p^{\vec{\varepsilon}}(y_i=\cdot
\;|z_{-n_1}^i),
$$
$$
\hat{x}_i^{\vec{\varepsilon}}=\hat{x}_i^{\vec{\varepsilon}}(\hat{z}_{-n_2}^i)=p^{\vec{\varepsilon}}(y_i=\cdot
\;|\hat{z}_{-n_2}^i).
$$
Then we have
$$
x_{i+1}^{\vec{\varepsilon}}=f_{z_{i+1}}^{\vec{\varepsilon}}(x_i^{\vec{\varepsilon}}),
\qquad
\hat{x}_{i+1}^{\vec{\varepsilon}}=f_{z_{i+1}}^{\vec{\varepsilon}}(\hat{x}_i^{\vec{\varepsilon}}).
$$
From (\ref{contract-1}) and (\ref{confine-orbit}),
%
it follows that there exists a positive constant $L$ independent of
$n_1$ and $n_2$ such that
\begin{equation}   \label{F-rho}
|x_0^{\vec{\varepsilon}}-\hat{x}_0^{\vec{\varepsilon}}| \leq L
\rho^n.
\end{equation}

Naturally
\begin{equation}
\label{amalg} p^{\vec{\varepsilon}}(z_0|z_{-n}^{-1}) = \sum_{\{y_0:
\Phi(y_0) = z_0\}} \sum_{y_{-1}}
\Delta^{\vec{\varepsilon}}(y_{-1},y_0)p^{\vec{\varepsilon}}(y_{-1}|z_{-n}^{-1}).
\end{equation}
Then, there is a positive constant $L'$, independent of $n_1, n_2$,
such that
 \begin{equation}   \label{F-rho'}
|p^{\vec{\varepsilon}}(z_0|z_{-n_1}^{-1})-
p^{\vec{\varepsilon}}(\hat{z}_0|\hat{z}_{-n_2}^{-1})| \leq L'
\rho^n.
\end{equation}
Since $\Delta^{\vec{\varepsilon}_0}$ satisfies conditions $1$ and
$2$, $p^{\vec{\varepsilon}}(z_0|z_{-n}^{-1})$ is bounded away from
$0$, uniformly in $\vec{\varepsilon} \in \Omega_{\mathbb{C}}$, $n$
and choices of $z_{-n}^{-1}$; thus there is a positive constant
$L''$, independent of $n_1, n_2$, such that
 \begin{equation}   \label{F-rho''}
|\log p^{\vec{\varepsilon}}(z_0|z_{-n_1}^{-1})- \log
p^{\vec{\varepsilon}}(\hat{z}_0|\hat{z}_{-n_2}^{-1})| \leq L''
\rho^n.
\end{equation}

Since for each $y \in \{1, \ldots, B\}$, $p^{\vec{\varepsilon}}(y)$
is analytic, from
$$
p^{\vec{\varepsilon}}(z)=\sum_{\Phi(y)=z} p^{\vec{\varepsilon}}(y),
$$
we deduce that $p^{\vec{\varepsilon}}(z)$ is analytic. Furthermore
since $ p^{\vec{\varepsilon}}(z_0|z_{-n}^{-1})$ is analytic on
$\Omega_{\mathbb{C}}$, we conclude $p^{\vec{\varepsilon}}(
z_{-n}^0)$ is analytic on $\Omega_{\mathbb{C}}$.

Choose $\sigma$ so that
$$
1 < \sigma < 1/\rho.
$$
If $r$ and $R$
are chosen sufficiently small,
then
\begin{equation}
\label{conditional-complex}
\sum_{z_0}| p^{\vec{\varepsilon}}(z_0|z_{-n}^{-1})| \leq \sigma,
~\eps \in \Omega_{\mathbb{C}}(r) \mbox{ and all sequences } z
\end{equation}
and
\begin{equation}
\label{stationary-complex}
 \sum_{z_0}| p^{\vec{\varepsilon}}(z_0)
|\leq \sigma, ~\eps \in \Omega_{\mathbb{C}}(r) .
\end{equation}
Then we have
$$
\sum_{z_{-n-1}^0}
|p^{\vec{\varepsilon}}(z_{-n-1}^0)|=\sum_{z_{-n-1}^0}
|p^{\vec{\varepsilon}}(z_{-n-1}^{-1})
p^{\vec{\varepsilon}}(z_0|z_{-n-1}^{-1})|  \leq \sum_{z_{-n-1}^{-1}}
|p^{\vec{\varepsilon}}(z_{-n-1}^{-1})| \sum_{z_0}
|p^{\vec{\varepsilon}}(z_0|z_{-n-1}^{-1})| \leq \sigma
\sum_{z_{-n}^0} |p^{\vec{\varepsilon}}(z_{-n}^{0})|,
$$
implying
\begin{equation} \label{bound-sigma} \sum_{z_{-n-1}^0}
|p^{\vec{\varepsilon}}(z_{-n-1}^0)|  \leq \sigma^{n+2}.
\end{equation}

Let
$$
H_n^{\vec{\varepsilon}}(Z)=-\sum_{z_{-n}^0}
p^{\vec{\varepsilon}}(z_{-n}^0) \log
p^{\vec{\varepsilon}}(z_0|z_{-n}^{-1})
$$
and
$$
\rho_1=\rho \delta < 1,
$$
then we have
$$
|H_{n+1}^{\vec{\varepsilon}}(Z)-H_n^{\vec{\varepsilon}}(Z)|=|\sum_{z_{-n-1}^0}
p^{\vec{\varepsilon}}(z_{-n-1}^0) \log
p^{\vec{\varepsilon}}(z_0|z_{-n-1}^{-1})-\sum_{z_{-n}^0}
p^{\vec{\varepsilon}}(z_{-n}^0) \log
p^{\vec{\varepsilon}}(z_0|z_{-n}^{-1})|
$$
$$
=|\sum_{z_{-n-1}^0} p^{\vec{\varepsilon}}(z_{-n-1}^0) (\log
p^{\vec{\varepsilon}}(z_0|z_{-n-1}^{-1})- \log
p^{\vec{\varepsilon}}(z_0|z_{-n}^{-1}))| \leq \sigma^2 L'' \rho_1^n;
$$
here the latter inequality follows from (\ref{F-rho''}) and
(\ref{bound-sigma}). Thus, for $m > n$,
$$
 |H_{m}^{\vec{\varepsilon}}(Z)-H_n^{\vec{\varepsilon}}(Z)| \leq
 \sigma^2 L''(\rho_1^n + \ldots +  \rho_1^{m-1}) \leq
 \frac{\sigma^2L''\rho_1^n}{1 - \rho_1}.
$$
This establishes the uniform convergence of
$H_n^{\vec{\varepsilon}}(Z)$ to a limit
$H_{\infty}^{\vec{\varepsilon}}(Z)$.  By Theorem $2.4.1$
of~\cite{ta02}, the uniform limit of complex analytic functions on a
fixed complex neighborhood is analytic on that neighborhood, and so
$H_{\infty}^{\vec{\varepsilon}}(Z)$ is analytic on
$\Omega_{\mathbb{C}}$.

For real $\vec{\varepsilon}$, $H_{\infty}^{\vec{\varepsilon}}(Z)$
coincides with the entropy rate function $H(Z^{\vec{\varepsilon}})$,
and so Theorem~\ref{main} follows.

\begin{exmp}  \label{bsc}
Consider a binary symmetric channel with crossover probability
$\varepsilon$. Let $\{Y_n\}$ be the input Markov chain with the
transition matrix
\begin{equation}   \label{transition}
\Pi=\left[\begin{array}{cc}
       \pi_{00}&\pi_{01}\\
       \pi_{10}&\pi_{11}\\
       \end{array}\right].
\end{equation}
At time $n$ the channel can be characterized by the following
equation
$$
Z_n=Y_n \oplus E_n,
$$
where $\oplus$ denotes binary addition, $E_n$ denotes the i.i.d.
binary noise with $p_E(0)=1-\varepsilon$ and $p_E(1)=\varepsilon$,
and $Z_n$ denotes the corrupted output. Then $(Y_n, E_n)$ is jointly
Markov, so $\{Z_n\}$ is a hidden Markov chain with the corresponding
$$
\Delta =
\left [ \begin{array}{cccc}
\pi_{00} (1-\varepsilon) & \pi_{00} \varepsilon & \pi_{01} (1-\varepsilon) & \pi_{01} \varepsilon\\
\pi_{00} (1-\varepsilon) & \pi_{00} \varepsilon & \pi_{01} (1-\varepsilon) & \pi_{01} \varepsilon\\
\pi_{10} (1-\varepsilon) & \pi_{10} \varepsilon & \pi_{11} (1-\varepsilon) & \pi_{11} \varepsilon\\
\pi_{10} (1-\varepsilon) & \pi_{10} \varepsilon & \pi_{11} (1-\varepsilon) & \pi_{11} \varepsilon\\
\end{array} \right ];
$$
here, $\Phi$ maps states $1$ and $4$ to $0$ and maps states $2$ and
$3$ to $1$. This class of hidden Markov chains has been studied
extensively (e.g.,~\cite{ja04},~\cite{or03}).

By Theorem~\ref{main}, when $\varepsilon$ and $\pi_{ij}$'s are
positive, the entropy rate $H(Z)$ is analytic as a function of
$\varepsilon$ and $\pi_{ij}$'s.   This still holds when $\varepsilon
= 0$ and the $\pi_{ij}$'s are positive, because in this case, we
have
$$
\Delta =
\left [ \begin{array}{cccc}
\pi_{00}  & 0 & \pi_{01}  & 0\\
\pi_{00} & 0 & \pi_{01}  & 0\\
\pi_{10}  & 0 & \pi_{11}  & 0\\
\pi_{10}  & 0 & \pi_{11}  & 0\\
\end{array} \right ].
$$

\end{exmp}

\section{Domain of Analyticity}   \label{domain}

Suppose $\Delta$ is analytically parameterized by a vector variable
$\vec{\varepsilon}$, and Conditions 1 and 2 in Theorem~\ref{main}
are satisfied at $\vec{\varepsilon}=\vec{\varepsilon}_0$. In
principle, the proof of Theorem~\ref{main} determines a neighborhood
$\Omega_{\mathbb{C}}(r)$ of $\vec{\varepsilon}_0$ on which the
entropy rate is analytic.  Specifically, if one can find $\rho, r$
and $R$ such that all of the following hold, then the entropy rate
is analytic on $\Omega_{\mathbb{C}}(r)$.

\begin{enumerate}
\item
Find $\rho$ such that each $f^{\eps_0}_a$ is a Euclidean
$\rho$-contraction on each $W_b$. Then choose positive $r, R$ such
that for all $\vec{\eps} \in \Omega_{\mathbb{C}}(r)$, each
$f_a^{\vec{\eps}}$ is a Euclidean $\rho$-contraction on each
$N_b(R)$ (see (\ref{contract-1})).
\item
Next find $r$ smaller (if necessary) such that for all
$\vec{\varepsilon} \in \Omega_{\mathbb{C}}(r)$, the image of the
stationary vector of $\Delta^{\vec{\eps}}$, under any composition of
the mappings $\{f^\eps_a\}$, stays within $\cup_b N_b(R)$ (see
(\ref{confine-orbit})).  Note that the argument in the proof shows
that this holds if (\ref{bounded-1}) and (\ref{bounded-2}) hold.
\item
Finally, find $r,R$ such that the sum of the absolute values of the
complexified conditional probabilities, conditioned on any given
past symbol sequence, is $< 1/\rho$, and similarly for the sum of
the absolute values of the complexified stationary probabilities
(see (\ref{conditional-complex}) and (\ref{stationary-complex})).
\end{enumerate}

In fact, the proof shows that one can always find such $\rho, r,R$,
but in condition $1$ above one may need to replace $f_a$'s by all
$n$-fold compositions of the $f_a$'s, for some $n$.

Recall from Example~\ref{bsc} the family of hidden Markov chains
$Z^\eps$ determined by passing a binary Markov chain through a
binary symmetric channel with cross-over probability $\varepsilon$.
Recall that $H(Z^\eps)$ is an analytic function of $\varepsilon$ at
$\varepsilon=0$ when the Markov transition probabilities are all
positive.  We shall determine a complex neighborhood of $0$ such
that the entropy rate, as a function of $\varepsilon$, is analytic
on this neighborhood.

Let $u_n=p(y_n=0|z_1^n)$ and $v_n=p(y_n=1|z_1^n)$. For $z_{n+1}=1$
we have
$$
u_{n+1}=\frac{\varepsilon(\pi_{00} u_n+\pi_{10}
v_n)}{\varepsilon(\pi_{00} u_n+\pi_{10} v_n)+ (1-\varepsilon)
(\pi_{01} u_n +\pi_{11} v_n)},
$$
$$
v_{n+1}=\frac{(1-\varepsilon)(\pi_{01} u_n+\pi_{11}
v_n)}{\varepsilon(\pi_{00} u_n+\pi_{10} v_n)+ (1-\varepsilon)
(\pi_{01} u_n +\pi_{11} v_n)}.
$$
Since $u_n+v_n=1$, $u_{n+1}$ is a function of $u_n$; let $g_1$
denote this function.

For $z_{n+1}=0$ we have
$$
u_{n+1}=\frac{(1-\varepsilon)(\pi_{00} u_n+\pi_{10}
v_n)}{(1-\varepsilon)(\pi_{00} u_n+\pi_{10} v_n)+ \varepsilon
(\pi_{01} u_n +\pi_{11} v_n)},
$$
$$
v_{n+1}=\frac{\varepsilon(\pi_{01} u_n+\pi_{11}
v_n)}{(1-\varepsilon)(\pi_{00} u_n+\pi_{10} v_n)+ \varepsilon
(\pi_{01} u_n +\pi_{11} v_n)}.
$$
Again, $u_{n+1}$ is a function of $u_n$; let $g_0$ denote this
function.

And for the conditional probability, we have
$$
p(z_n=0|z_1^{n-1}) = ((1-\varepsilon) \pi_{00} + \varepsilon
\pi_{01}) u_n+((1-\varepsilon) \pi_{10} + \varepsilon \pi_{11}) v_n.
$$
Since $u_n+v_n=1$, $p(z_n=0|z_1^{n-1})$ is a function of $u_n$; let
$r_0$ denote this function. And
$$
p(z_n=1|z_1^{n-1}) =(\varepsilon \pi_{00} + (1-\varepsilon)
\pi_{01}) u_n+(\varepsilon \pi_{10} + (1-\varepsilon) \pi_{11}) v_n.
$$
Again, $p(z_n=1|z_1^{n-1})$ is a function of $u_n$; let $r_1$ denote
this function.

Note that $g_0, g_1, r_0, r_1$ are all implicitly parameterized by
$\varepsilon$. The stationary vector $(\pi_0, \pi_1)$ of $Y$, which
doesn't depend on $\varepsilon$, is equal to
$(\pi_{10}/(\pi_{10}+\pi_{01}), \pi_{01}/(\pi_{10}+\pi_{01}))$.

We shall choose $\rho$ with $0 < \rho < 1$, $r > 0$ and $R > 0$ such
that for all $\varepsilon$ with $|\varepsilon| < r$
\begin{enumerate}
\item  $g_0$ and $g_1$ are $\rho$-contraction mappings
on $R$-neighborhoods of $0$ and $1$ in the complex plane,
\item
the set of all $\{g_{a_n} \circ g_{a_{n-1}} \circ \cdots \circ
g_{a_1}(\pi_0) \}$) are within the $R$-neighborhoods of $0$ and $1$,
\item and $|r_0(u)|+|r_1(u)| < 1/\rho$ for $u$ in $R$-neighborhoods of $0$ and $1$ in the complex plane.
\end{enumerate}
By the general principle above, the entropy rate should be analytic
on $|\varepsilon| < r$.

More concretely, condition $1$, $2$ and $3$ translate to (here $\rho
< 1$):
\begin{enumerate}
\item $|g'_0(u)| < \rho$, $|g'_1(u)| < \rho$ on
($|\varepsilon| < r$ and $|u| < R$) and ($|\varepsilon| < r$ and
$|1-u| < R$),
\item $\max{\{|g_0(0)-1|, |g_0(1)-1|, |g_1(0)|, |g_1(1)|\}} < R
(1-\rho)$ on $|\varepsilon| < r$ (this follows from
(\ref{bounded-1}); (\ref{bounded-2}) is trivial since the stationary
vector of $Y$ doesn't depend on $\varepsilon$),
\item $|r_0(u)|+|r_1(u)| < 1/\rho$ on ($|\varepsilon| < r$ and $|u| < R$) and
($|\varepsilon| < r$ and $|1-u| < R$).
\end{enumerate}

A straightforward computation shows that the following conditions
guarantee conditions $1$, $2$, $3$:
$$
0 < \frac{\sqrt{r (|-\pi_{00} \pi_{11}+\pi_{10} \pi_{11}+\pi_{10}
\pi_{01} -\pi_{10} \pi_{11}| r+|(\pi_{00} \pi_{11}+\pi_{10}
\pi_{01})|)}}{\pi_{11}-|\pi_{10}-\pi_{11}|
r-(|\pi_{00}-\pi_{10}-\pi_{01}+\pi_{11}| r + |\pi_{01}-\pi_{11}|) R}
< \sqrt{\rho},
$$

$$
0 < \frac{\sqrt{r (|-\pi_{00} \pi_{11}+\pi_{10} \pi_{11}+\pi_{10}
\pi_{01} -\pi_{10} \pi_{11}| r+|(\pi_{00} \pi_{11}+\pi_{10}
\pi_{01})|)}}{\pi_{01}-|\pi_{00}-\pi_{01}|
r-(|\pi_{00}-\pi_{10}-\pi_{01}+\pi_{11}| r + |\pi_{01}-\pi_{11}|) R}
< \sqrt{\rho},
$$

$$
0 < \frac{\sqrt{r(|-\pi_{11} \pi_{00} + \pi_{01} \pi_{00} + \pi_{01}
\pi_{10} -\pi_{01} \pi_{00}| r + |\pi_{11} \pi_{00} -\pi_{01}
\pi_{10}|)}}{\pi_{00}-|\pi_{01}-\pi_{00}|
r-(|\pi_{00}-\pi_{10}+\pi_{11}-\pi_{01}| r + |\pi_{10}-\pi_{00}|) R}
< \sqrt{\rho},
$$

$$
0 < \frac{\sqrt{r(|-\pi_{11} \pi_{00} + \pi_{01} \pi_{00} + \pi_{01}
\pi_{10} -\pi_{01} \pi_{00}| r + |\pi_{11} \pi_{00} -\pi_{01}
\pi_{10}|)}}{\pi_{10}-|\pi_{11}-\pi_{10}|
r-(|\pi_{00}-\pi_{10}+\pi_{11}-\pi_{01}| r + |\pi_{10}-\pi_{00}|) R}
< \sqrt{\rho},
$$

$$
0 < \frac{r \pi_{00}}{\pi_{01}-|\pi_{00}-\pi_{01}|r} < R (1-\rho), 0
< \frac{r \pi_{10}}{\pi_{11}-|\pi_{10}-\pi_{11}|r} < R (1-\rho),
$$

$$
0 < \frac{r \pi_{11}}{\pi_{10}- |\pi_{11}-\pi_{10}|r} < R (1-\rho),
0 < \frac{r \pi_{01}}{\pi_{00} -|\pi_{01}-\pi_{00}|r} < R (1-\rho),
$$

$$
(|\pi_{00}-\pi_{01}-\pi_{10}+\pi_{11}| r + |\pi_{01}-\pi_{11}|) R +
|\pi_{10}-\pi_{11}| r+ \pi_{11},
$$
$$
+(|\pi_{01}-\pi_{00}+\pi_{10}-\pi_{11}| r + |\pi_{00}-\pi_{10}|) R +
|\pi_{11}-\pi_{10}| r + \pi_{10} < 1/\rho,
$$

$$
(|\pi_{10}-\pi_{11}-\pi_{00}+\pi_{01}| r + |\pi_{11}-\pi_{01}|) R +
|\pi_{00}-\pi_{01}| r+ \pi_{01}
$$
$$
+(|\pi_{11}-\pi_{10}+\pi_{00}-\pi_{01}| r + |\pi_{10}-\pi_{00}|) R +
|\pi_{01}-\pi_{00}| r + \pi_{00} < 1/\rho.
$$

In other words, for given $\rho$ with $0 < \rho < 1$, choose $r$ and
$R$ to satisfy all the constraints above. Then the entropy rate is
an analytic function of $\varepsilon$ on $|\varepsilon| < r$.

\section{Relaxed Conditions}
\label{relaxed}

We do not know a complete set of necessary and sufficient conditions
on $\Delta$ and $\Phi$ that guarantee analyticity of entropy rate.
However, in this section, we show how the hypotheses in
Theorem~\ref{main} can be relaxed and still guarantee analyticity.
We then give several examples. In Section~\ref{boundary}, we do give
a a complete set of necessary and sufficient conditions for a very
special class of hidden Markov chains.

In this section, we assume that $\Delta$ has a simple maximum
eigenvalue $1$; this implies that $\Delta$ has a unique stationary
vector $\vec{s}$.

For a mapping $f$ from $W_b$ to $W$ and $w \in W_b$. Let $f'$ denote
the first derivative of $f$ at $w$ restricted to the subspace
spanned by directions parallel to the simplex $W_b$ and let $\|
\cdot \|$ denote the Euclidean norm of a linear mapping. We say that
$\{f_a: a \in A\}$ is \emph{eventually contracting} at $w \in W_b$
if there exists $n$ such that for any $a_0, a_1, \cdots, a_n \in A$,
$\|(f_{a_n} \circ f_{a_{n-1}} \circ \cdots \circ f_{a_0})'(w)\|$ is
strictly less than $1$. We say that $\{f_a: a \in A\}$ is \emph{
contracting} at $w \in W_b$ if it is \emph{eventually contracting}
at $w$ with $n=0$.   Using the mean value theorem, one can show that
if  $\{f_a: a \in A\}$ is \emph{ contracting}  at each $w$ in a
compact convex subset $K$ of $W_b$ then each $f_a$ is a contraction
mapping on $K$.

Let $L$ denote the limit set of $\{(f_{a_n} \circ f_{a_{n-1}} \circ
\cdots \circ f_{a_0})(\vec{s})\}$.

\begin{thm} \label{general}
If at $\Delta=\hat\Delta$,
\begin{enumerate}
\item $1$ is a simple eigenvalue for $\hat\Delta$,
\item For every $a$ and all $w$ in $L$, $r_a(w) > 0$,
\item For every $b$, $\{f_a: a \in A\}$ is eventually contracting at
all $w$ in the convex hull of the intersection of $L$ and $W_b$,
\end{enumerate}
then $H(Z)$ is analytic at $\Delta=\hat\Delta$.
\end{thm}

\begin{proof}

Let $\mathcal{X}$ denote the right infinite shift space
$\{a_0^{\infty}: a_i \in A\}$. Let $L_{\delta}$ be the set of all
points in $W$ of distance at most $\delta$ from $L$.  Choose
$\delta$ so small that
\begin{itemize}
\item For every $a \in A$ and $w$ in $L_{\delta}$, $r_a(w) > 0$ --
and --
\item For every $b$, $\{f_a: a \in A\}$ is eventually contracting at
all $w$ in the convex hull of the intersection of $L_{\delta}$ and
$W_b$.
\end{itemize}

Since the convex hull $K_\delta$ of the intersection of $L_{\delta}$
and $W_b$ is compact, there exists $n$ such that for any $a_0, a_1,
\cdots, a_n \in A$ and any $w \in K_\delta$, $\|(f_{a_n} \circ
f_{a_{n-1}} \circ \cdots \circ f_{a_0})'(w)\|$ is strictly less than
$1$. For simplicity, we may assume that for each $a$, $\{f_a\}$ is
contracting on $K_\delta$, and so each $f_a$ is a contraction
mapping on $K_\delta$.   Since $L_\delta \subseteq K_\delta$, it
follows that $f_a(L_\delta) \subseteq L_\delta$ and so each $f_a$ is
a contraction mapping on $L_\delta$.

For any $c_0^{\infty} \in \mathcal{X}$, there exists $n$ such that
$\{(f_{c_n} \circ f_{c_{n-1}} \circ \cdots \circ f_{c_0})(\vec{s})\}
\in L_{\delta}$. Let $\mathcal{X}_{c_0^{\infty}}^n$ denote the
cylinder set $\{a_0^{\infty}: a_0=c_0, a_1=c_1, \cdots, a_n=c_n\}$.
Since $\{f_a: a \in A\}$ is a contraction mapping on $L_\delta$, we
conclude that for any $a_0^{\infty} \in
\mathcal{X}_{c_0^{\infty}}^n$ and all $m \geq n$, $\{(f_{a_m} \circ
f_{a_{m-1}} \circ \cdots \circ f_{a_0})(\vec{s})\} \in L_{\delta}$.
By the compactness of $\mathcal{X}$, we can find finitely many such
cylinder sets to cover $\mathcal{X}$. Consequently we can find $n$
such that for any $a_0^{\infty} \in \mathcal{X}$ and any $m \ge n$ ,
we have $\{(f_{a_m} \circ f_{a_{m-1}} \circ \cdots \circ
f_{a_0})(\vec{s})\} \in L_{\delta}$. We can now apply the proof of
Theorem~\ref{main} -- namely, we can use the contraction (along any
symbolic sequence $z_{-n}^0$) to extend $H_n(Z)=H(Z_0|Z^{-1}_{-n})$
from real to complex and prove the uniform convergence of $H_n(Z)$
to $H(Z)$ in complex parameter space.
\end{proof}

\begin{rem}
\label{perron}~

(1) If $\hat \Delta$ has a strictly positive column (or more
generally, there is a $j$ such that for all $i$, there exists $n$
such that $\hat \Delta^n_{ij} >0$),  then condition $1$ of
Theorem~\ref{general} holds by Perron-Frobenius theory.

(2) If for each symbol $a$, $\hat \Delta_a$ is row allowable (i.e.,
no row is all zero), then $r_a(w)> 0$ for all $w \in W$ and so
condition $2$ of Theorem~\ref{general} holds.
\end{rem}
\bigskip

Theorem~\ref{general} relaxes the positivity assumptions of
Theorem~\ref{main}. Indeed given conditions $1$ and $2$ of
Theorem~\ref{main}, by Remark~\ref{perron}, conditions 1 and 2 of
Theorem~\ref{general} hold. For condition $3$ of
Theorem~\ref{general}, first observe that $L$ is contained in
$\cup_b f_b(W)$. Using the equivalence of the Euclidean metric and
the Hilbert metric, Proposition~\ref{eventually-contract} shows that
for every $b$, $\{f_a: a \in A\}$ is eventually contracting on
$f_b(W)$, which is a convex set containing the intersection of $L$
and $W_b$.

Theorem~\ref{general} also applies to many cases not covered by
Theorem~\ref{main}.  For instance, suppose that some column of
$\hat \Delta$ is strictly positive and each $\hat \Delta_a$ is row
allowable. By Remark~\ref{perron}, Theorem~\ref{general} applies
whenever we can guarantee condition 3.  For this, it is sufficient
to check that for each $a, b$, $f_a$ is a contraction, with respect
to the Euclidean metric, on the convex hull of the intersection of
$L$ with each $W_b$. This can be done by explicitly computing
derivatives.

\begin{exmp}
Consider a hidden Markov chain $Z$ defined by :
$$
\hat \Delta =  \left[ \begin{array}{cccc}
a_{11} &a_{12} & a_{13}& a_{14}\\
a_{21} &a_{22} & a_{23}& a_{24}\\
a_{31} &a_{32} & a_{33}& a_{34}\\
a_{41} &a_{42} & a_{43}& a_{44}\\
\end{array} \right],
$$
with $\Phi(1)=\Phi(2)=0$ and $\Phi(3)=\Phi(4)=1$. We assume that
some column of $\hat \Delta$ is strictly positive and both $\hat
\Delta_0$ and $\hat \Delta_1$  are row allowable.

Parameterize $W_0$ by $(y, 1-y, 0, 0)$ and parameterize $W_1$ by
$(0, 0, y, 1-y)$ (with  $y \in [0,1]$). We can explicitly compute
the derivatives of $f_0$ and $f_1$ with respect to $y$:
$$
f'_0|_{(y, 1-y, 0, 0)}=\frac{a_{11} a_{22}-a_{12}
a_{21}}{((a_{11}+a_{12}-a_{21}-a_{22})y+a_{21}+a_{22})^2},
$$
$$
f'_0|_{(0, 0, y, 1-y)}=\frac{a_{31} a_{42}-a_{32}
a_{41}}{((a_{31}+a_{32}-a_{41}-a_{42})y+a_{41}+a_{42})^2},
$$
$$
f'_1|_{(y, 1-y, 0, 0)}=\frac{a_{13} a_{24}-a_{14}
a_{23}}{((a_{13}+a_{14}-a_{23}-a_{24})y+a_{23}+a_{24})^2},
$$
$$
f'_1|_{(0, 0, y, 1-y)}=\frac{a_{33} a_{44}-a_{34}
a_{43}}{((a_{33}+a_{34}-a_{43}-a_{44})y+a_{43}+a_{44})^2},
$$
Note that the row allowability condition guarantees that the
denominators in these expressions never vanish.

Choose $a_{ij}$'s such that each of these derivatives is less than
1; then we conclude that the entropy rate is analytic at
$\hat\Delta$. One way to do this is to make each of the $2 \times 2$
upper/lower left/right matrices singular.

Or choose the $a_{ij}$'s such that
$$ \hat \Delta = \left [ \begin{array}{cccc}
\alpha_1& * & \beta_1 & 0\\
0       & \alpha_2& 0 & \beta_2\\
\lambda_1& * & \eta_1 & 0\\
0       & \lambda_2& 0 & \eta_2\\
\end{array} \right ]
$$
where $0< \alpha_1 < \alpha_2$, $0< \beta_1 < \beta_2$, $0 <
\lambda_1 < \lambda_2$, $0 < \eta_1 < \eta_2$ and $*$ denote a real
positive number. Let $(s_2, s_4)$ be the Perron eigenvalue of the
stochastic matrix:
$$
\left [ \begin{array}{cc}
\alpha_2 & \beta_2\\
\lambda_2 & \eta_2\\
\end{array} \right ].
$$
Then $\vec{s}=(0, s_2, 0, s_4)$ is the stationary vector of $\Delta$
corresponding to the simple eigenvalue $1$. Let $w_0=(0, 1, 0, 0)$
and $w_1=(0, 0, 0, 1)$. One checks that for $n \geq 0$, $f_{a_n}
\circ f_{a_{n-1}} \circ \cdots \circ f_{a_0}(\vec{s})=w_{a_n}$.
Therefore $L$ consists of $\{w_0, w_1\}$. Using the expressions
above, we see that $$ f'_0|_{w_0}=\alpha_1/\alpha_2 < 1,
f'_0|_{w_1}=\lambda_1/\lambda_2 <1,
$$
$$
f'_1|_{w_0}=\beta_1/\beta_2 <1, f'_1|_{w_1}=\eta_1/\eta_2 < 1.
$$
So, $f_0$ and $f_1$ are contraction mappings at $\{w_0, w_1\}$, and
so condition 3 holds. Thus, the entropy rate $H(Z)$ is analytic at
$\hat\Delta$.
\end{exmp}

\section{Hidden Markov Chains with Unambiguous Symbol}
\label{boundary}

\begin{de}
A symbol $a$ is called {\em unambiguous} if $\Phi^{-1}(a)$ contains
only one element.
\end{de}

\begin{rem}
Note that unambiguous symbol is referred to as ``singleton clump''
in some ergodic theory work, such as~\cite{pe03}.
\end{rem}

When an unambiguous symbol is present, the entropy rate can be
expressed in a simple way: letting $a_1$ be an unambiguous symbol,
\begin{equation}
\label{ent-unambig-0} H(Z)=\sum_{a_{i_j} \neq a_1} p(a_{i_n}
a_{i_{n-1}} \cdots a_{i_2} a_1) H(z|a_{i_n} a_{i_{n-1}} \cdots
a_{i_2} a_1).
\end{equation}

In this section, we focus on the case of a binary hidden Markov
chain, in which 0 is unambiguous.  Then, we can rewrite
(\ref{ent-unambig-0}) as
\begin{equation} \label{ent-unambig}
H(Z^{\vec\eps}) = p^{\vec\eps}(0) H^{\vec\eps}(z|0)+
p^{\vec\eps}(10) H^{\vec\eps}(z|10)+ \cdots + p^{\vec\eps}(1^{(n)}0)
H^{\vec\eps}(z|1^{(n)}0)+ \cdots ,
\end{equation}
where $1^{(n)}$ denotes the sequence of $n$ 1's and
$$
H^{\vec\eps}(z|1^{(n)}0) = - p^{\vec\eps}(0|1^{(n)}0) \log
p^{\vec\eps}(0|1^{(n)}0) - p^{\vec\eps}(1|1^{(n)}0) \log
p^{\vec\eps}(1|1^{(n)}0).
$$

\begin{exmp}  \label{exmp1}
Fix $a, b, \ldots, h >0$ and for $\eps \ge 0$ let
$$
\Delta(\eps)=\left [ \begin{array}{ccc}
                \eps&a - \eps &b\\
                g&c&d\\
                h&e&f\\
                \end{array} \right ].
$$
Assume $a, b, \ldots, h >0$ are chosen such that $\Delta(\eps)$ is
stochastic. The symbols of the Markov chain are the matrix indices
$\{1,2,3\}$. Let $Z^\eps$ be the binary hidden Markov chain defined
by: $\Phi(1)=0$ and $\Phi(2)=\Phi(3)=1$. We claim that $H(Z^\eps)$
is not analytic at $\eps = 0$.

Let $\pi(\eps)$ be the stationary vector of $\Delta(\eps)$ (which is
unique since $\Delta(\eps)$ is irreducible). Observe that
$$
p^\eps(0) = \pi_1(\eps),~ p^\eps(00)=\pi_1 (\eps) \varepsilon,
$$
and for $n \geq 1$.
$$
p^\eps(1^{(n)}0)=\pi_1(\eps) (a - \eps, b) {\left[\begin{array}{cc}
                                              c&d\\
                                              e&f\\
                                       \end{array}\right]}^{n-1}
                                       \left( \begin{array}{c}
                                              1\\
                                              1\\
                                       \end{array} \right).
$$
Since $\Delta(\eps)$ is irreducible, $\pi(\eps)$ is analytic in
$\eps$ and positive. Now,
\begin{equation}
\label{z0} p^\eps(0) H^\eps(z|0) = - p^\eps(00) \log p^\eps(0|0) -
p^\eps(10) \log p^\eps(1|0).
\end{equation}
The first term in (\ref{z0}) is
$$
- p^\eps(00) \log p^\eps (0|0)= - \pi_1 (\eps) \varepsilon \log
\varepsilon,
$$
which is not analytic (or even differentiable at $\eps =0$). The
second term in (\ref{z0}) is
$$
- p^\eps (10) \log p^\eps (1|0)= - \pi_1(\eps) (a -\eps + b) \log
(\pi_1(\eps) (a -\eps + b)),
$$
which is analytic at $\eps =0$.  Thus, $H^\eps(z|0)$ is not analytic
at $\eps = 0$.  Similarly it can be shown that all of the terms of
(\ref{ent-unambig}), other than $H^\eps(z|0)$, are analytic at $\eps
= 0$. Since the matrix
$$
\left[\begin{array}{cc}
                                              c&d\\
                                              e&f\\
                                       \end{array}\right]
$$
has spectral radius $ < 1$, the terms of (\ref{ent-unambig}) decay
exponentially; it follows that the infinite sum of these terms is
analytic. Thus,  $H(Z^\eps)$ is the sum of two functions of $\eps$,
one of which is analytic and the other is not analytic at $\eps =0$.
Thus, $H(Z^\eps)$ is not analytic at $\eps = 0$.

\end{exmp}

\begin{exmp}  \label{exmp2}
Fix $a, b, \cdots, g > 0$ and consider the stochastic matrix
$$
\Delta(\varepsilon)=\left [ \begin{array}{ccc}
                e&a&b\\
                f-\varepsilon&c&\varepsilon\\
                g&0&d\\
                \end{array} \right ].
$$
The symbols of the Markov chain are the matrix indices $\{1, 2,
3\}$. Again let $Z^{\varepsilon}$ be the binary hidden Markov chain
defined by  $\Phi(1)=0$ and $\Phi(2)=\Phi(3)=1$. We show that
$H(Z^{\varepsilon})$ is analytic at $\varepsilon=0$ when $c \neq d$,
and not analytic when $c = d$. Note that
$$
p^{\varepsilon}(0)=\pi_1 (\varepsilon),
$$
and for $n \geq 1$.
$$
p^{\varepsilon}(1^{(n)}0) =\pi_1 (\varepsilon) (a, b)
{\left[\begin{array}{cc}
                                              c&\varepsilon\\
                                              0&d\\
                                       \end{array}\right]}^{n-1}
                                       \left( \begin{array}{c}
                                              1\\
                                              1\\
                                       \end{array} \right).
$$
When $c \neq d$, we assume $c > d$, then
$$
{\left[\begin{array}{cc}
                                              c&\varepsilon\\
                                              0&d\\
                                       \end{array}\right]}^n=\left[\begin{array}{cc}
                                              c^n&\varepsilon c^{n-1} \frac{1-(d/c)^{n}}{1-d/c}\\
                                              0&d^n\\
                                       \end{array}\right].
$$
Since $\Delta(\varepsilon)$ is irreducible, $\pi(\varepsilon)$ is
analytic in $\varepsilon$ and positive. Simple computation leads to:
$$
p^{\varepsilon}(1|1^{(n)}0)= (a c^n+ a \varepsilon c^{n-1}
\frac{1-(d/c)^{n}}{1-d/c} + b d^n)/(a c^{n-1}+ a \varepsilon c^{n-2}
\frac{1-(d/c)^{n-1}}{1-d/c}+ b d^{n-1} )
$$
$$
=(a c^2+ a \varepsilon c \frac{1-(d/c)^{n}}{1-d/c} + bd^2
(d/c)^{n-2})/(ac+\varepsilon \frac{1-(d/c)^{n-1}}{1-d/c} + bd
(d/c)^{n-2}),
$$
and
$$
p^{\varepsilon}(0|1^{(n)}0)=((f-\eps) a c^{n-1}+ g (a\varepsilon
c^{n-2} \frac{1-(d/c)^{n-1}}{1-d/c} + b d^{n-1}))/(a c^{n-1}+a
\varepsilon c^{n-2} \frac{1-(d/c)^{n-1}}{1-d/c}+ b d^{n-1} )
$$
$$
=((f -\eps) a c + g(a \varepsilon \frac{1-(d/c)^{n-1}}{1-d/c} + b d
(d/c)^{n-2}) )/(ac + \varepsilon \frac{1-(d/c)^{n-1}}{1-d/c} + bd
(d/c)^{n-2}).
$$
In this case all terms are analytic. Again since
$$
\left [ \begin{array}{cc}
              c&\varepsilon\\
              0&d\\
\end{array} \right ]
$$
has spectral radius $< 1$, the term $p^{\varepsilon}(1^{(n)}0)
H^{\varepsilon}(z|1^{(n)}0)$ is exponentially decaying with respect
to $n$. Therefore the infinite sum of these terms is also analytic,
and so the entropy rate is a real analytic function of
$\varepsilon$.

When $c=d$, we have
$$
p^{\varepsilon}(1|1^{(n)}0)= (ac^{n+1}+a \varepsilon (n+1) c^n + b
c^{n+1})/(a c^n + a \varepsilon n c^{n-1}+ b c^n)
$$
$$
=(ac^2+ a \varepsilon (n+1)c+ bc^2)/( ac + a \varepsilon n +bc),
$$
and
$$
p^{\varepsilon}(0|1^{(n)}0)=((f -\eps) a c^n+ga \varepsilon n
c^{n-1}+ g b c^n)/(a c^n+ a \varepsilon n c^{n-1} + b c^n)
$$
$$
=((f -\eps) a c+ ga \varepsilon n + gbc)/(ac + a \varepsilon n +
bc).
$$
For any $n$, consider a small neighborhood $N_n$ of $-(a+b)c/an$ in
$\mathbb{C}$ such that $-(a+b)c/a j \in N_n$ only holds for $j=n$.
When $\varepsilon \to -(a+b)c/an$, the complexified term
$p^{\varepsilon}(1^{(n)}0) H^{\varepsilon}(z|1^{(n)}0) \to \infty$.
Meanwhile, the sum of all the other terms can be analytically
extended to $N_n$ (from any path $I$ from a positive $\varepsilon$
to $-(a+b)c/an$ with $-(a+b)c/aj \notin I$ for $j \neq n$). Thus, by
the uniqueness of analytic continuation of $H(Z^{\varepsilon})$, we
conclude that $H(Z^{\varepsilon})$ blows up when one approaches
$-(a+b)c/an$ and therefore is not analytic at $\varepsilon=0$
(although it is smooth from the right at $\varepsilon=0$).
\end{exmp}

The two examples above show that under certain conditions the
entropy rate of a binary hidden Markov chain with unambiguous symbol
can fail to be analytic at the boundary.  We now show that these
examples typify all the types of failures of analyticity at the
boundary (in the case of a binary hidden Markov chains with an
unambiguous symbol).

We will need the following result.

\begin{lem}  \label{Senata}
Let $A(\vec\eps)$ be an analytic parameterization of complex
matrices. Let $\lambda$ be the spectral radius of $A(\vec\eps_0)$.
Then for any $\eta > 0$, there exists a complex neighborhood
$\Omega$ of $\vec\eps_0$ and positive constant $C$ such that for all
$\vec\eps \in \Omega$ and all $i,j,k$
$$
|A_{ij}^k(\vec\varepsilon)| \leq C (\lambda +\eta)^k.
$$
\end{lem}

\begin{proof}
Following ~\cite{se80}, we consider
$$
(I- z A)^{-1}=I + zA + z^2 A^2+ \cdots.
$$
And
$$
(I-z
A)^{-1}=\frac{Adj(I-zA)}{\det(I-zA)}=\frac{Adj(I-zA)}{(1-\lambda_1
z)(1-\lambda_2 z) \cdots (1-\lambda_n z)},
$$
where $\lambda_1, \ldots, \lambda_n$ are the eigenvalues of $A$. So
every entry of $(I-z A)^{-1}$ takes the form:
$$
(p_0+p_1 z+ \cdots + p_m z^m) \prod_{j=1}^n \sum_{i=0}^{\infty}
\lambda_j^i z^i
$$
$$
=\sum_{k=0}^\infty \sum_{u=0}^m p_u \sum_{i_1+i_2+\cdots+i_n=k-u}
\lambda_1^{i_1} \lambda_2^{i_2} \cdots \lambda_n^{i_n} z^k.
$$
Since the eigenvalues of a complex matrix vary continuously with
entries, the lemma follows.
\end{proof}

Now let $S(n)$ denote the set of all the $n \times n$ complex
matrices with isolated (in modulus) maximum eigenvalue.
\begin{lem}  \label{connected}
$S(n)$ is connected.
\end{lem}

\begin{proof}
let $A, B \in S(n)$, then we consider their Jordan forms:
$$
A= U \diag(\lambda_1, C) U^{-1}, \qquad B= V \diag(\eta_1, D)
V^{-1},
$$
here $\lambda_1, \eta_1$ are maximum eigenvalues for $A, B$,
respectively, $C, D$ correspond to other Jordan blocks, and $U, V
\in GL(n, \mathbb{C})$ (here $GL(n, \mathbb{C})$ denotes the set of
all the $n \times n$ nonsingular complex matrices). Since $GL(n,
\mathbb{C})$ is connected~\cite{on93}, it suffices to prove that
there is a path in $S(n)$ from $\diag(\lambda_1, C)$ to
$\diag(\eta_1, D)$. This is straightforward: first connect
$\diag(\lambda_1, C)$ to $\diag(\eta_1, \eta_1/\lambda_1 C)$ by a
continuous rescaling; then connect $\eta_1/\lambda_1 C$ to $D$ by
the path $t\eta_1/\lambda_1 C + (1-t)D$ (the path
$\diag(\eta_1,t\eta_1/\lambda_1 C +(1-t)D )$ stays within $S(n)$
since the matrices along this path are upper triangular with all
diagonal entries, except $\eta_1$, of modulus less than $|\eta_1|$).
\end{proof}

For a complex analytic function $f(z_1, z_2, \cdots, z_n)$, let
$V(f)$ denote the ``hypersurface'' defined by $f$, namely
$$
V(f)=\{(z_1, z_2, \cdots, z_n) \in \mathbb{C}^n: f(z_1, z_2, \cdots,
z_n)=0\}.
$$
Now let $\Omega$ denote a connected open set in $\mathbb{C}^n$. It
is well known that the following Lemma holds (for completeness, we
include a brief proof).
\begin{lem}   \label{nzf}
$\Omega \backslash V(f)$ is connected.
\end{lem}

\begin{proof}
For simplicity, we first assume $\Omega$ is a ball $B_r(z_0)$ (here
$z_0 \in \mathbb{C}^n$ is the center of the ball and $r$ is the
radius, i.e., $B_r(z_0)=\{z \in \mathbb{C}^n: |z-z_0| < r\}$) in
$\mathbb{C}^n$. For any two distinct point $P, Q \in \Omega
\backslash V(f)$, consider the ``complex line''
$$
L_{\mathbb{C}}^{PQ}=\{zP+(1-z)Q: z \in \mathbb{C}\}.
$$
$L_{\mathbb{C}}^{PQ} \cap V(f) \cap \Omega$ consists of only
isolated points (A non-constant one variable complex analytic
function must have isolated zeros in the complex plane~\cite{sh92}).
It then follows that for the compact real line segment:
$$
L_{\mathbb{R}}^{PQ}=\{tP+(1-t)Q: t \in [0,1]\},
$$
$L_{\mathbb{R}}^{PQ} \cap V(f) \cap \Omega$ consists of only
finitely many points. Certainly one can choose an arc in
$L_{\mathbb{C}}^{PQ} \cap \Omega$ to avoid these points and connect
$P$ and $Q$. This implies that $\Omega \backslash V(f)$ is
connected.

In the general case, $\Omega$ is a connected open set in
$\mathbb{C}^n$. Let $I$ be an arc in $\Omega$ connecting $P$ and
$Q$, and let $\{B_{r_j}(z_j)\}$ be a collection of balls covering
$I$ such that each $B_{r_j}(z_j) \cap B_{r_{j+1}}(z_{j+1}) \ne
\phi$. Pick a point $P_j$ in $B_{r_j}(z_j) \cap
B_{r_{j+1}}(z_{j+1})$ such that $P_j \in \Omega\backslash V(f)$.
Applying the same argument as above to every ball $B_{r_j}(z_j)$, we
see that $P$ is connected to $Q$ in  $\Omega\backslash V(f)$ through
the points $P_j$'s. Thus we prove the lemma.
\end{proof}

\begin{thm}
\label{unambig-bdry} Let $\Delta$ be an irreducible stochastic $d
\times d$ matrix.  Write $\Delta$ in the form:
\begin{equation}
\label{form-2} \Delta =\left [ \begin{array}{cc}
                a &r \\
                c &B\\
       \end{array} \right ]
\end{equation}
where $a$ is a scalar and $B$ is a $(d-1) \times (d-1)$ matrix. Let
$\Phi$ be the function defined by $\Phi(1) = 0$, and $\Phi (2) =
\cdots =\Phi(d) = 1$. Then for any parametrization
$\Delta(\vec\eps)$ such that $\Delta(\vec\eps_0) = \Delta$, letting
$Z^{\vec\eps}$ denote the hidden Markov chain defined by
$\Delta(\vec\eps)$ and $\Phi$, $H(Z^{\vec\eps})$ is analytic at
$\vec\varepsilon_0$ if and only if
\begin{enumerate}
\item $a > 0$, and $r B ^jc  > 0$ for $j=0, 1, \cdots$. \item The
maximum eigenvalue of $B $ is simple and strictly greater in
absolute value than the other eigenvalues of $B$.
\end{enumerate}
\end{thm}

\begin{proof}
\medskip

\textbf{Proof of sufficiency.}

We write
\begin{equation}
\label{form} \Delta(\vec\eps)=\left [ \begin{array}{cc}
                a(\vec\eps) &r(\vec\eps) \\
                c(\vec\eps) &B (\vec\eps)\\
       \end{array} \right ],
\end{equation}
where $a(\vec\eps)$ is a scalar and $B(\vec\eps)$ is a $(d-1) \times
(d-1)$ matrix.

Since $\Delta(\vec\eps_0)$ is stochastic and irreducible, its
spectral radius is 1, and  1 is a simple eigenvalue of $\Delta$.
Thus, if $\Omega$ is sufficiently small, for all $\vec\eps \in
\Omega$, any fixed row $\pi(\vec\eps)=(\pi_1(\vec\varepsilon),
\pi_2(\vec\varepsilon), \cdots, \pi_d(\vec\varepsilon))$ of $Adj(I -
\Delta(\vec\eps))$ is a left eigenvector of $\Delta(\vec\eps)$
associated with eigenvalue $1$ and is an analytic function of
$\vec\eps$. Normalizing, we can assume that $\pi(\vec\eps)
\mathbf{1} = \mathbf{1}$, $\pi(\vec\eps)$ is analytic in $\vec\eps$,
and $\pi(\vec\eps_0)
>0$.

The entries of $r(\vec\eps), B(\vec\eps),$  and $c(\vec\eps)$ are
real analytic in $\vec\eps$ and can be extended to complex analytic
functions in a complex neighborhood $\Omega$ of $\vec\eps_0$. Thus,
for all $n$, $\pi_1(\vec\eps)r(\vec\eps)B(\vec\eps)^{n-1}
\mathbf{1}$ and $\pi_1(\vec\eps)r(\vec\eps)B(\vec\eps)^{n-1}
c(\vec\eps)$ can be extended to complex analytic functions on
$\Omega$ (in fact, each of these functions is a polynomial in
$\vec\eps$).

Since $B(\vec\eps_0)$ is a proper sub-matrix of the irreducible
stochastic matrix $\Delta(\vec\eps_0)$, its spectral radius is
strictly less than 1. Thus, by Lemma~\ref{Senata}, there exists $0 <
\lambda^* < 1$ and a constant $C_1
>0$, such that for some complex neighborhood $\Omega$ of $\vec\eps_0$, all
$\vec\eps \in \Omega$, and all n,
$$
|B^n_{ij}(\vec\eps)|  < C_1(\lambda^*)^n.
$$
Since $\pi_1(\vec\eps)$,  $r(\vec\eps)$ and $c(\vec\eps)$ are
continuous in $\vec\eps$, there is a constant $C_2 >0$ such that for
all $\vec\eps \in \Omega$ and all $n$:
\begin{equation}
\label{decay} |\pi_1(\vec\eps)r(\vec\eps)B(\vec\eps)^{n} \mathbf{1}|
< C_2 (\lambda^*)^n.
\end{equation}

We will need the following result, proven in
Appendix~\ref{lem-bounds}.

\begin{lem}
\label{bounds} Let
$$
a(\vec\eps,n) \equiv \frac{\pi_1(\vec\eps)r(\vec\eps)B(\vec\eps)^n
\mathbf{1}}{\pi_1(\vec\eps)r(\vec\eps)B(\vec\eps)^{n-1} \mathbf{1}}
$$
and
$$
b(\vec\eps,n) \equiv
\frac{\pi_1(\vec\eps)r(\vec\eps)B(\vec\eps)^{n-1}
c(\vec\eps)}{\pi_1(\vec\eps)r(\vec\eps)B(\vec\eps)^{n-1}
\mathbf{1}}.
$$
For a sufficiently small neighborhood  $\Omega$ of $\vec\eps_0$,
both $a(\vec\eps,n)$ and $b(\vec\eps,n)$ are bounded from above and
away from zero, uniformly in $\vec\eps \in \Omega$ and $n$.
\end{lem}

Define
$$
H^{\vec\eps}_n  = -  a(\vec\eps,n) \log a(\vec\eps, n) - b(\vec\eps,
n) \log b(\vec\eps,n),
$$
where $a(\vec\eps, n)$ and $b(\vec\eps, n)$ are as in
Lemma~\ref{bounds}. Choosing $\Omega$ to be a smaller neighborhood
of $\vec\eps_0$, if necessary, $a(\vec\eps, n)$ and $b(\vec\eps, n)$
are constrained to lie in a closed disk not containing $0$. Thus for
all $n$, $H^{\vec\eps}_n$ is an analytic function of $\vec\eps$,
with $|H^{\vec\eps}_n |$ bounded uniformly in $\vec\eps \in \Omega$
and $n$. Since $\pi_1(\vec\eps)r(\vec\eps)B(\vec\eps)^{n-1} {\bf 1}$
is analytic on $\Omega$ and exponentially decaying (by
(\ref{decay})), the infinite series
\begin{equation}   \label{HZ}
H^{\vec\eps}(Z) = \pi_1(\vec\eps) H^{\vec\eps}_0 +
\pi_1(\vec\eps)r(\vec\eps) {\bf 1} H^{\vec\eps}_1 + \cdots +
\pi_1(\vec\eps)r(\vec\eps)B(\vec\eps)^{n-1} {\bf 1} H^{\vec\eps}_n +
\cdots
\end{equation}
converges uniformly on $\Omega$ and thus defines an analytic
function on $\Omega$.

Note that for $\vec\eps \ge 0$,
\begin{equation}
\label{ppi1} p^{\vec\eps} (1^{(n)} 0) =
\pi_1(\vec\eps)r(\vec\eps)B(\vec\eps)^{n-1} {\bf 1}
\end{equation}
and
\begin{equation}
\label{ppi0} p^{\vec\eps} (01^{(n)}0) = \pi_1(\vec\eps)r(\vec\eps
)B(\vec\eps)^{n-1} c(\vec\eps).
\end{equation}
By (\ref{ppi1}), (\ref{ppi0}), and the expression for entropy rate
in the case of an unambiguous symbol (given at the beginning of this
section), $H^{\vec\eps}(Z)$ agrees with the entropy rate when
$\Delta(\vec\eps) \ge 0$, as desired.
\medskip

\begin{rem}
We show how sufficiency relates to Theorem~\ref{general}. Namely,
the assumptions in Theorem~\ref{unambig-bdry} imply those of
Theorem~\ref{general}. Condition $1$ of Theorem~\ref{general}
follows from the fact that $\Delta$ is assumed irreducible. For
conditions 2 and 3 of Theorem~\ref{general}, one first notes that
the image of $f_0$ is a single point $W_0$, and the $f_1$-orbit of
$W_0$ and $f_1$-orbit of $\vec{s}$ converge to a point $p_1$. It
follows that $L$ is the union of $W_0$, the $f_1$-orbit of $W_0$ and
$p_1$. The assumptions in Theorem~\ref{unambig-bdry}. imply that
$r_a
>0$ on $L$ (i.e., condition 2 of Theorem~\ref{general} holds) and that for
sufficiently large $n$, the $n$-fold composition of $f_1$ is
contracting on the convex hull of the intersection of $L$ and $W_1$
(so condition 3 of Theorem~\ref{general} holds). To see the latter,
one uses the ideas in the proof of sufficiency.
\end{rem}

\textbf{Proof of necessity}

We first consider condition $2$. We shall use the natural
parameterization and view $H(Z)$ as a function of $\Delta$, or more
precisely of $(B, r)$. Note that there is a one-to-one
correspondence between $\Delta$ and $(B, r)$; we shall use this
correspondence throughout the proof.

Suppose $\Delta$ doesn't satisfy condition $2$, however $H(Z)$ is
analytic at $\Delta$ with respect to the natural parameterization.
In other words, suppose there exists a complex neighborhood
$N_{\Delta}$ of $\Delta$ (here $N_{\Delta}$ corresponds to $N_B
\times N_r$ where $N_B$ is neighborhood of $B$ and $N_r$ is
neighborhood of $r$) such that $H(Z)$ can be analytically extended
to $N_{\Delta}$, while the corresponding $B$ doesn't have isolated
(in modulus) maximum eigenvalue.

We first claim there exists $\tilde{\Delta} \in N_{\Delta}$ with
$\tilde{r} \tilde{B}^k \mathbf{1}=0$, here $\tilde{r}$ and
$\tilde{B}$ correspond to $\tilde{\Delta}$ and $\tilde{B}$ has
distinct eigenvalues (in modulus). Indeed we can first (for
simplicity) perturb $\Delta$ to $\tilde{\Delta}$ such that the
corresponding $\tilde{B}$ has distinct eigenvalues in modulus. Then
$$
\tilde{B}=\tilde{U} \diag(\tilde{\lambda}_1, \tilde{\lambda}_2,
\cdots, \tilde{\lambda}_{d-1}) \tilde{U}^{-1}
$$
$$
=(\tilde{v}_1, \tilde{v}_2, \cdots, \tilde{v}_{d-1})
\diag(\tilde{\lambda}_1, \tilde{\lambda}_2, \cdots,
\tilde{\lambda}_{d-1}) (\tilde{w}_1^t, \tilde{w}_2^t, \cdots,
\tilde{w}_{d-1}^t)^t
$$
where $|\tilde{\lambda}_1| > |\tilde{\lambda}_2| > \cdots >
|\tilde{\lambda}_{d-1}|$, and $\tilde{v}_i, \tilde{w}_i$'s are
appropriately scaled right and left eigenvectors of $\tilde{B}$,
respectively. Then we have
$$
r \tilde{B}^k \mathbf{1}= r \tilde{v}_1 \tilde{w}_1 \mathbf{1}
\tilde{\lambda}_1^k+ r \tilde{v}_2 \tilde{w}_2 \mathbf{1}
\tilde{\lambda}_2^k+ \cdots+ r \tilde{v}_{d-1} \tilde{w}_{d-1}
\mathbf{1} \tilde{\lambda}_{d-1}^k.
$$
Further consider a perturbation of $B$ from
$$
\tilde{B}=\tilde{U} \diag(\tilde{\lambda}_1, \tilde{\lambda}_2,
\cdots, \tilde{\lambda}_{d-1}) \tilde{U}^{-1}
$$
to
$$
\tilde{B}=V\tilde{U} \diag(\tilde{\lambda}_1, \tilde{\lambda}_2,
\cdots, \tilde{\lambda}_{d-1}) \tilde{U}^{-1}V^{-1},
$$
where $V$ is a complex matrix close to the $(d-1) \times (d-1)$
identity matrix $I_{d-1}$. So we can pick $V$ such that $\tilde{v}_1
\tilde{w}_1 V^{-1} \mathbf{1} \neq 0$, $\tilde{v}_1 \tilde{w}_1
V^{-1} \tilde{c} \neq 0$, $\tilde{v}_2 \tilde{w}_2 V^{-1} \mathbf{1}
\neq 0$. Clearly $\tilde{v}_1 \tilde{w}_1 V^{-1} \mathbf{1}$ is not
proportional to $\tilde{v}_2 \tilde{w}_2 V^{-1} \mathbf{1}$. Then by
a further perturbation of $r$ to $\tilde{r}$, we can simultaneously
require that $\tilde{r} \tilde{v}_1 \tilde{w}_1 \mathbf{1} \neq 0$,
$\tilde{r} \tilde{v}_1 \tilde{w}_1 \tilde{c} \neq 0$, $\tilde{r}
\tilde{v}_2 \tilde{w}_2 \mathbf{1} \neq 0$, $ |\tilde{r} \tilde{v}_1
\tilde{w}_1 \mathbf{1}| \neq |\tilde{r} \tilde{v}_2 \tilde{w}_2
\mathbf{1}|$, where we redefine $\tilde{v}_i=V \tilde{v}_i$ and
$\tilde{w}_i= \tilde{w}_i V^{-1}$. For any $\theta$ and $\eta > 0$,
it can be checked that
$$
\bigcup_{k=0}^{\infty} \{z^k: |z- e^{i \theta}| < \eta \} =
\mathbb{C} \backslash \{0\}.
$$
Since $\tilde{\lambda}_2$ is a perturbation of $\tilde{\lambda}_1$,
it follows that for large enough $k$, one can perturb
$\tilde{\lambda}_2$ to satisfy the equation
$$
{\left(\tilde{\lambda}_2/ \tilde{\lambda}_1
\right)}^k=\frac{-\tilde{r} \tilde{v}_1 \tilde{w}_1
\mathbf{1}-\tilde{r} \tilde{v}_3 \tilde{w}_3 \mathbf{1}
(\tilde{\lambda}_3/\tilde{\lambda}_1)^k - \cdots-\tilde{r}
\tilde{v}_{d-1} \tilde{w}_{d-1} \mathbf{1}
(\tilde{\lambda}_{d-1}/\tilde{\lambda}_1)^k }{\tilde{r} \tilde{v}_2
\tilde{w}_2 \mathbf{1}},
$$
with $|\tilde{\lambda}_2| \neq |\tilde{\lambda}_1|$ and
$|\tilde{\lambda}_2|$ strictly greater than $|\tilde{\lambda}_j|$
for $j \geq 3$. Thus we prove the claim.

We now pick a positive matrix $\hat{\Delta} \in N_{\Delta}$ with
corresponding $\hat{r}$ and $\hat{B}$. We then pick $\tilde{\Delta}
\in N_{\Delta}$ with corresponding $\tilde{r}$ and $\tilde{B}$ (with
distinct eigenvalues in modulus) such that $\tilde{r}
\tilde{B}^{k_1} \mathbf{1}=0$ for some $k_1$, and we can further
require that $\tilde{r} \tilde{v}_1 \tilde{w}_1 \mathbf{1} \neq 0$,
$\tilde{r} \tilde{v}_1 \tilde{w}_1 \tilde{c} \neq 0$ (see the proof
for the previous claim), where as before, $\tilde{v}_1, \tilde{w}_1$
are eigenvectors corresponding to the largest eigenvalue of
$\tilde{B}$. According to Lemma~\ref{connected}, there is an arc
$I_1 \subset S(d-1)$ connecting $\hat{B}$ to $\tilde{B}$; we then
connect $\hat{r}$ and $\tilde{r}$ using an arc $I_2$ in
$\mathbb{C}^{d-1}$. According to Lemma~\ref{nzf}, we can choose the
arc $I=(I_1, I_2)$ to avoid the hypersurface $V((r v_1 w_1
\mathbf{1})(r v_1 w_1 c)) \subset \mathbb{C}^{(d-1)^2} \times
\mathbb{C}^{d-1}$; in other words, we can assume that along the path
$I$, $r v_1 w_1 \mathbf{1} \neq 0$ and $r v_1 w_1 c \neq 0$; here
$v_1, w_1, c$ are determined by the variable matrix $B$ along the
path $I_1$ and $r$ is the variable point along path $I_2$ (we remind
the reader that the coordinates of $v_1$ and $w_1$ are all analytic
functions of the entries of $B$). We then claim that there is a
neighborhood $N_I$ of $I$ such that $V_k \cap N_I \neq \phi$ and
$W_k \cap N_I \neq \phi$ hold for only finitely many $k$, where
$V_k=\{(B, r): r B^k \mathbf{1} =0\}$ and $W_k=\{(B, r): r B^k c
=0\}$. Indeed for any $\Delta \in I$ with corresponding $B \in
S(d-1)$, by the Jordan form we have
$$
r B^k \mathbf{1}=r v_1 w_1 \mathbf{1} \lambda_1^k +o(\lambda_1^k),
$$
where $\lambda_1$ is the isolated maximum eigenvalue and $v_1, w_1$
are appropriately scaled right and left eigenvectors of $B$,
respectively. Since $r v_1 w_1 \mathbf{1} \neq 0$ on $I$, there
exists a complex connected neighborhood $N_I$ of $I$ such that $r
v_1 w_1 \mathbf{1} \neq 0$ on $N_I$ and $r v_1 w_1 \mathbf{1}
\lambda_1^k$ dominates uniformly on $N_I$ (see Lemma~\ref{Senata}).
Consequently, $ | r B^k \mathbf{1} | > 0$ on $N_I$ for large enough
$k$. In other words, $V_k \cap N_I \neq \phi$ holds for only
finitely many $k$. Similarly since $r v_1 w_1 c \neq 0$ on $I$,
there exists a complex neighborhood $N_I$ of $I$ (here we use the
same notation for a possibly different neighborhood) such that $W_k
\cap N_I \neq \phi$ holds only for finitely many $k$. From now on,
we assume such $k$'s are less than some $K$, which depends on $N_I$.

We claim that we can further choose $I$ and find a new neighborhood
$N_I$ in $\mathbb{C}^{d-1} \times S(d-1)$ of $I$ such that $V_k \cap
N_I \neq \phi$ holds only for $k=k_1$ and $W_k \cap N_I =\phi$ for
all $k$. Consider $\tilde{\Delta}$ with corresponding $\tilde{B}$,
let $F_i=F_i(\tilde{B})= \{r: r \tilde{B}^i \mathbf{1}=0\}$, which
is a hyperplane orthogonal to the vector $\tilde{B}^i \mathbf{1}$ in
$\mathbb{C}^{d-1}$. Similarly we define $G_i=G_i(\tilde{B})= \{r: r
\tilde{B}^i \tilde{c}=0\}$. Recall that $\tilde{B}=\tilde{U}
\diag(\lambda_1, \lambda_2, \cdots, \lambda_{d-1}) \tilde{U}^{-1}$;
we can require that $\tilde{U}^{-1} \mathbf{1}$ has no zero
coordinates by a small perturbation of $\tilde{U}$ if necessary. We
then show that $F_i$'s and $G_j$'s define different hyperplanes in
$\mathbb{C}^{d-1}$. Indeed suppose $F_i=F_j$. It follows that
$\tilde{U} \diag(\tilde{\lambda}_1^i, \tilde{\lambda}_2^i, \cdots,
\tilde{\lambda}_{d-1}^i) \tilde{U}^{-1} \mathbf{1}$ is proportional
to $\tilde{U} \diag(\tilde{\lambda}_1^j, \tilde{\lambda}_2^j,
\cdots, \tilde{\lambda}_{d-1}^j) \tilde{U}^{-1} \mathbf{1}$. It then
follows that $(\tilde{\lambda}_1^i, \tilde{\lambda}_2^i, \cdots,
\tilde{\lambda}_{d-1}^i)$ is proportional to $(\tilde{\lambda}_1^j,
\tilde{\lambda}_2^j, \cdots, \tilde{\lambda}_{d-1}^j)$. However
since not all eigenvalues have the same modulus, this implies that
$i=j$. With a perturbation of $\tilde{c}$ (equivalently a
perturbation of row sums of $\tilde{B}$), if necessary, we conclude
that the $F_i$'s and $G_i$'s determine different hyperplanes, i.e.,
$ F_i \neq F_j$, $G_i \neq G_j$ for $i \neq j \leq K$, and $F_i \neq
G_j$ for all $i, j$. Thus, with a perturbation of $\tilde{r}$ if
necessary, we can choose a new $\tilde{\Delta}$ contained in
$V_{k_1}$, but not contained in any $V_k$ with $k \neq k_1$ or $W_k$
for all $k$. Again by Lemma~\ref{nzf}, one can choose a new $I$
inside original $N_I$, connecting $\hat{\Delta}$ and
$\tilde{\Delta}$, to avoid all $V_k$'s and $W_k$'s except $V_{k_1}$,
then choose a smaller new neighborhood $N_I$ of the new $I$ to make
sure that $V_k \cap N_I \neq \phi$ only holds for $k=k_1$ and $W_k
\cap N_I = \phi$ for all $k$.

Since the perturbed complex matrix $B$ still has spectral radius
strictly less than $1$, all the complexified terms in the entropy
rate formula (see (\ref{HZ})) with $k \neq k_1$ are exponentially
decaying and thus sum up to an analytic function  on $N_I$.
(i.e., the sum of these terms can be analytically continued to
$N_I$), while the unique analytic extension of the $k_1$-th term on
$N_I$ blows up as one approaches $V_{k_1} \cap N_I$ from
$\hat{\Delta}$. Again by the uniqueness of analytic extension of
$H(Z)$ on $N_I$, this would be a contradiction to the assumption
that $H(Z)$ is analytic at $\Delta$ (here we are applying the
uniqueness theorem of analytic continuation of a function of several
complex variables, see page $21$ in~\cite{sh92}). Thus we prove the
necessity of condition $2$.

We now consider condition $1$. Suppose $\Delta$ doesn't satisfies
condition $1$, namely $a=0$ or $r B^k c = 0$ for some $k$, however
$H(Z)$ is analytic at $\Delta$. With the proof above for the
necessity of condition $2$, we can now assume the corresponding $B
\in S(d-1)$.

If $a=0$, consider any perturbation of $\Delta$ to $\Delta_1$ such
that $\tilde{B} \in S(d-1)$, $\tilde{r} \tilde{v}_1 \tilde{w}_1
\mathbf{1} \neq 0$, $\tilde{r} \tilde{v}_1 \tilde{w}_1 \tilde{c}
\neq 0$, $\tilde{r} \tilde{B}^k \mathbf{1} \neq 0$ and $\tilde{r}
\tilde{B}^k \tilde{c} \neq 0$ for all $k$ (here we follow the
notation as in the proof of necessity of condition $2$). Then using
similar arguments, we can prove the sum of all the terms except the
first term in the entropy rate formula (see (\ref{HZ})) can be
analytically extended to $\tilde{\Delta}$. However this implies that
$a \log a$ is a well-defined analytic function on some neighborhood
of $0$ in $\mathbb{C}$, which is a contradiction. Similar arguments
can be applied to the case that $r B^k c=0$ for some $k$'s. Thus we
prove the necessity of condition $1$.
\end{proof}

\section{Analyticity of a Hidden Markov Chain in a Strong Sense} \label{xxx-II}

In this section, we show that if $\Delta$ is analytically
parameterized by a real variable vector $\vec{\varepsilon}$, and at
$\vec{\varepsilon}_0$, $\Delta$ satisfies conditions $1$ and $2$ of
Theorem~\ref{main}, then the hidden Markov chain {\em itself} is a
real analytic function of $\vec{\varepsilon}$ at
$\vec{\varepsilon}_0$ in a strong sense. We assume (for this section
only) that the reader is familiar with the basics of measure theory
and functional analysis~\cite{mu78, yo74, na81}. Our approach uses a
connection between the entropy rate of a hidden Markov chain and
symbolic dynamics explored in ~\cite{ma84a1}.

Let $\mathcal{X}$ denote the set of left infinite sequences with
finite alphabet. A cylinder set is a set of the form:
$(\{x_{-\infty}^0: x_0=z_0, \cdots, x_{-n}=z_{-n}\})$. The Borel
sigma-algebra is the smallest sigma-algebra containing the cylinder
sets. A Borel probability measure (BPM) $\nu$ on $\mathcal{X}$ is a
measure on the Borel measurable sets of $\mathcal{X}$ such that
$\nu(\mathcal{X})=1$. Such a measure is uniquely determined by its
values on the cylinder sets.

For real $\vec{\varepsilon}$, consider the measure
$\nu^{\vec{\varepsilon}}$ on $\mathcal{X}$ defined by:
\begin{equation}
\label{ppp} \nu^{\vec{\varepsilon}}(\{x_{-\infty}^0: x_0=z_0,
\cdots, x_{-n}=z_{-n}\})=p^{\vec{\varepsilon}}(z_{-n}^0).
\end{equation}

Note that $H(Z)$ can be rewritten as
\begin{equation}
\label{ent-form} H^{\vec{\varepsilon}}(Z)=\int -\log
p^{\vec{\varepsilon}}(z_0|z_{-\infty}^{-1})
d\nu^{\vec{\varepsilon}}.
\end{equation}

Usually, the Borel sigma-algebra is defined to be the smallest
sigma-algebra containing the open sets; in this case, the open sets
are defined by the metric: for any two elements $\xi$ and $\eta$ in
$\mathcal{X}$, define $d(\xi, \eta)=2^{-k}$ where $k=\inf\{|i|:
\xi_i \neq \eta_i \}$. The metric space $(\mathcal{X}, d)$ is
compact.

Let $C(\mathcal{X})$ be the space of real-valued continuous
functions on $\mathcal{X}$.  Then $C(\mathcal{X})$  is a Banach
space (i.e., complete normed linear space) with the sup norm
$||f||_\infty = \sup\{|f(x)|: x \in   \mathcal{X}\}$.  Then any BPM
$\nu$ acts as a bounded linear functional on $C(\mathcal{X})$,
namely $\nu(f) =\int f d\nu$.  As such, the set of BPM's is a subset
of the dual space, $C(\mathcal{X})^*$, which is itself a Banach
space; the norm of a BPM $\nu$ is defined: $||\nu||= \sup_{\{f \in
C(\mathcal{X}): ||f||_\infty =1\}} \int f d\nu$. In fact, since
$\mathcal{X}$ is compact,  $C(\mathcal{X})^*$ is the linear span of
the BPM's.

It makes sense to ask if $\vec{\varepsilon} \mapsto
\nu^{\vec{\varepsilon}}$ is analytic as a mapping from the parameter
space to $C(\mathcal{X})^*$; by definition, this would mean that
$\nu^{\vec{\varepsilon}}$ can be expressed as a power series in the
coordinates of $\vec{\varepsilon}$.  However, as the following
example shows, this mapping is not even continuous.

Let $\mathcal{X}$ be the set of binary left infinite sequences. Let
$\nu_p$ denote the i.i.d. $(p,1-p)$ measure, with $0 < p <1$. Let
$$
S_p = \{ x \in X: \lim_{n \to \infty} (1/n)( \log p_{x_1} + \ldots +
\log p_{x_{-n}} ) = -p\log p -(1-p)\log (1-p)\}.
$$

Note that $S_p$ is a Borel measurable set. By the strong law of
large numbers,  $\nu_p(S_p) = 1$.  Clearly, for distinct $p$, $S_p$
are disjoint. Thus, for $q \ne p$,  $\nu_q(S_p) = 0$.

Any Borel measurable set $S$ can be approximated by a finite union
of cylinder sets in the following sense: given $\delta > 0$ and $p
\in (0,1)$, there is a finite union $C$ of cylinder sets such that
$|\nu_q(S) - \nu_q(C)| < \delta$ for all $q$ in a neighborhood of
$p$.  Applying this fact to $S= S_p$, and denoting $C_{(p,\delta)} =
C$, we obtain
$$
1 = \nu_p(S_p) - \nu_q(S_p) \leq  |
\nu_p(S_p)-\nu_p(C_{(p,\delta)})| + |\nu_p(C_{(p,\delta)}) -
\nu_q(C_{(p,\delta)})| + |\nu_q(C_{(p,\delta)}) - \nu_q(S_p)|
$$
$$
\leq 2 \delta +  |\nu_p(C_{(p,\delta)}) - \nu_q(C_{(p,\delta)})|.
$$
If $\delta < 1/2$, then  $\nu_q(C_{(p,\delta)})$ cannot converge to
$\nu_p(C_{(p,\delta)})$ as $q \to p$.   Since the characteristic
function of a finite union of cylinder sets is continuous, this
shows that the map $p \mapsto \nu_p$ from $\mathcal{R}$ to $C(X)^*$
is discontinuous.

On the other hand, using the work of Ruelle~\cite{ru78}, we now show
that $\vec{\varepsilon} \mapsto \nu^{\vec{\varepsilon}}$ is analytic
as a mapping from the parameter space to another natural space.

For $f \in C(\mathcal{X})$, define $ var_n (f)=\sup
\{|f(\xi)-f(\xi')| : \xi_{-i}=\xi'_{-i} \,\, \mbox{for} \,\, i \leq
n\}$. We denote by $F^{\theta}$ the subset of $f \in C(\mathcal{X})$
such that
$$
{\|f\|}_{\theta} \equiv \sup_{n \geq 0} (\theta^{-n} var_n (f)) < +
\infty.
$$
$F^{\theta}$ is a Banach space with the norm
$\|f\|=\max({|f|}_{\infty}, {\|f\|}_{\theta})$. Using complex
functions instead of real functions, one defines
$F_{\mathbb{C}}^{\theta}$ similarly.

In the following theorem, we prove the analyticity of a hidden
Markov chain in a strong sense.
\begin{thm}  \label{main-1}
Suppose that the entries of $\Delta$ are analytically parameterized
by a real variable vector $\vec{\varepsilon}$. If at
$\vec{\varepsilon}=\vec{\varepsilon}_0$, $\Delta$ satisfies
conditions $1$ and $2$ in Theorem~\ref{main}, then the mapping
$\vec{\varepsilon} \mapsto \log
p^{\vec{\varepsilon}}(z_0|z_{-\infty}^{-1})$ is analytic at
$\vec{\varepsilon}_0$ from the real parameter space to $F^{\rho}$
(here $\rho$ is the contraction constant in the proof of
Theorem~\ref{main}). Moreover the mapping $\vec{\varepsilon} \mapsto
\nu^{\vec{\varepsilon}}$ is analytic at $\vec{\varepsilon}_0$ from
the real parameter space to $(F^{\rho})^*$, the dual space (i.e.,
bounded linear functionals) on $F^\rho$.
\end{thm}

\begin{proof}
For complex $\vec{\varepsilon}$, by (\ref{F-rho''}), one shows that
$\log p^{\vec{\varepsilon}}(z_0|z_{-\infty}^{-1})$ can be defined on
$\Omega_{\mathbb{C}}$ as the uniform (in $ \vec{\varepsilon}$ and $z
\in \mathcal{X}$) limit of $\log
p^{\vec{\varepsilon}}(z_0|z_{-n}^{-1})$ as $n \to \infty$, and $\log
p(z_0|z_{-\infty}^{-1})$ belongs to $F^{\rho}$. By (\ref{x-i-0}),
(\ref{iter}), (\ref{init}) and (\ref{amalg}) it follows that $
p^{\vec{\varepsilon}}(z_0|z_{-n}^{-1})$  is analytic on
$\Omega_{\mathbb{C}}$. As a result of (\ref{F-rho''}), if $\Delta$
satisfies conditions $1$ and $2$, for fixed $z\in \mathcal{X}$,
$\log p^{\vec{\varepsilon}}(z_0|z_{-\infty}^{-1})$ is the uniform
limit of analytic functions and hence is analytic on
$\Omega_{\mathbb{C}}$ (see Theorem $2.4.1$ of~\cite{ta02}).

Using (\ref{F-rho''}) and the Cauchy integral formula in several
variables~\cite{ta02} (which expresses the derivative of an analytic
function at a point as an integral of a closed curve around the
point), we obtain the following.  There is a positive constant $C'$
such that whenever $z_{-\infty}^0 \stackrel{n}{\sim}
\hat{z}_{-\infty}^0$, for all $\vec{\varepsilon} \in
\Omega_{\mathbb{C}}$
\begin{equation}
\label{izabella} |D_{\vec{\varepsilon}}(\log
p^{\vec{\varepsilon}}(z_0|z_{-n_1}^{-1}))-D_{\vec{\varepsilon}}(\log
p^{\vec{\varepsilon}}(\hat{z}_0|\hat{z}_{-n_2}^{-1}))| \leq C'
\rho^n.
\end{equation}
Therefore for arbitrary yet fixed $ z_{-\infty}^0$, the components
of the derivatives of $\log
p^{\vec{\varepsilon}}(z_0|z^{-1}_{-\infty})$ with respect to
$\vec{\varepsilon}$ are also in $F_{\mathbb{C}}^{\rho}$.

Furthermore, we prove that the mapping $\vec{\varepsilon} \mapsto
\log p^{\vec{\varepsilon}}(z_0|z_{-\infty}^{-1})$ is complex
differentiable (therefore analytic) from $\Omega_{\mathbb{C}}$ to
$F_{\mathbb{C}}^{\theta}$. Let $f(\vec{\varepsilon}; \cdot)=\log
p^{\vec{\varepsilon}}(\cdot)$. It suffices to prove that
\begin{equation}  \label{differentiable-infty}
{\|f(\vec{\varepsilon}+\vec{h};\cdot)-f(\vec{\varepsilon};\cdot)-D_{\vec{\varepsilon}}
f|_{\vec{\varepsilon}}(\vec{h}; \cdot)\|}_{\infty} \leq o(\vec{h}).
\end{equation}
and
\begin{equation}  \label{differentiable}
{\|f(\vec{\varepsilon}+\vec{h};\cdot)-f(\vec{\varepsilon};\cdot)-D_{\vec{\varepsilon}}
f|_{\vec{\varepsilon}}(\vec{h}; \cdot)\|}_{\theta} \leq o(\vec{h}).
\end{equation}

Again applying the Cauchy integral formula in several variables, it
follows that there exists a positive constant $C''$ such that for
all $\vec{\varepsilon} \in \Omega_{\mathbb{C}}$ we have
\begin{equation}
\label{cauchy-upper}  |D^2_{\vec{\varepsilon}}
f|_{\vec{\varepsilon}}(\vec{h}, \vec{h}; z) | \leq C'' |\vec{h}\\|^2
\end{equation}
and whenever $z_{-\infty}^0 \stackrel{n}{\sim} \hat{z}_{-\infty}^0$,
\begin{equation}
\label{cauchy-rho} \int_0^1 (1-t) |(D_{\vec{\varepsilon}}^2
f|_{\vec{\varepsilon}}(\vec{h}, \vec{h}; z)-D^2_{\vec{\varepsilon}}
f|_{\vec{\varepsilon}}(\vec{h}, \vec{h}; \hat{z}))| dt \leq C''
|\vec{h}\\|^2 \rho^n,
\end{equation}
From the Taylor formula with integral remainder, we have:
 \begin{equation}
\label{Taylor} f(\vec{\varepsilon}+\vec{h}; z)-f(\vec{\varepsilon};
z)-D_{\vec{\varepsilon}} f|_{\vec{\varepsilon}}(\vec{h}; z
)=\int_0^1 (1-t) D_{\vec{\varepsilon}}^2
f|_{\vec{\varepsilon}+t\vec{h}}(\vec{h}, \vec{h}; z) dt.
\end{equation}
 To prove (\ref{differentiable-infty}), use
(\ref{cauchy-upper}) and (\ref{Taylor}).  To prove
(\ref{differentiable}), use (\ref{cauchy-rho}) and (\ref{Taylor}).
Therefore $ \vec{\varepsilon} \mapsto \log
p^{\vec{\varepsilon}}(\cdot)$ is  analytic as a mapping from
$\Omega_{\mathbb{C}}$ to $F^{\rho}_{\mathbb{C}}$. Restricting the
mapping $\vec{\varepsilon} \mapsto \log
p^{\vec{\varepsilon}}(z_0|z_{-\infty}^{-1})$ to the real parameter
space, we conclude that it is real analytic (as a mapping into
$F^{\rho}$). Using this and the theory of equilibrium
states~\cite{ru78}), the ``Moreover'' is proven in
Appendix~\ref{equilibrium}.
\end{proof}

\begin{co}
Suppose that at $\vec{\varepsilon}_0$, $\Delta$ satisfies conditions
$1$ and $2$ in Theorem~\ref{main}, and $\vec{\varepsilon} \mapsto
f^{\vec{\varepsilon}} \in F^{\rho}$ be analytic at
$\vec{\varepsilon}_0$, then $\vec{\varepsilon} \mapsto
\nu^{\vec{\varepsilon}} (f^{\vec{\varepsilon}})$ is analytic at
$\vec{\varepsilon}_0$. In particular, we recover Theorem~\ref{main}:
$\vec{\varepsilon} \mapsto H^{\vec{\varepsilon}}(Z)$ is analytic at
$\vec{\varepsilon}_0$.
\end{co}

\begin{proof}
The map
$$
\Omega \rightarrow F^\rho \times (F^\rho)^* \rightarrow \mathbb{R}
$$
$$
\vec{\varepsilon} \mapsto ( f^{\vec{\varepsilon}},
\nu^{\vec{\varepsilon}}) \mapsto
\nu^{\vec{\varepsilon}}(f^{\vec{\varepsilon}})
$$
is analytic at $\vec{\varepsilon}_0$, as desired.
\end{proof}

\vspace{1cm} \textbf{Acknowledgements:} We are grateful to Wael
Bahsoun, Joel Feldman, Robert Israel, Izabella Laba,   Erik
Ordentlich, Yuval Peres, Gadiel Seroussi, Wojciech Szpankowski and
Tsachy Weissman for helpful discussions.

\section*{Appendices}\appendix

\section{Proof of Proposition~\ref{Metric-Equivalence}}
\label{Metric-Equivalence1}
\begin{proof}
Without loss of generality, we assume $S$ is convex (otherwise
consider the convex hull of $S$). It follows from standard arguments
that max norm and sum norm are equivalent. More specifically, for
another metric $d_1$ defined by
$$
d_1(u, v)=\sqrt{\sum_{i \neq j \leq k}  \log^2 \left(
\frac{u_i/u_j}{v_i/v_j}. \right)},
$$
we have $d_\textbf{B} \sim d_1$. For metric $d_2$ defined by
$$
d_2(u, v)= \sqrt{\sum_{i \neq j \leq k} (u_i/u_j - v_i/v_j)^2}.
$$
Applying mean value theorem to $\log$ function, one concludes that $
d_1 \sim d_2$. Note that
\begin{eqnarray*}
u_i-v_i&=&\frac{u_i}{u_1+u_2+\cdots+u_k}-\frac{v_i}{v_1+v_2+\cdots+v_k}\\
       &=&\frac{1}{u_1/u_i+u_2/u_i+\cdots+u_k/u_i}-\frac{1}{v_1/v_i+v_2/v_i+\cdots+v_k/v_i}\\
\end{eqnarray*}
Applying the mean value theorem to function $f$, defined as
$$
f(x_1, x_2, \cdots, x_B)=\frac{1}{x_1+x_2+\cdots+x_k},
$$
we conclude that there exists $\xi \in S$ such that
$$
u_i-v_i=\nabla f|_{\xi} \cdot (u_1/u_i-v_1/v_i, \cdots,
u_k/u_i-v_k/v_i).
$$
It follows from Cauchy inequality that there exists a positive
constant $D_1$ such that
$$
d_{\mathbf{E}}(u, v) < D_1 d_2(u, v).
$$
Similarly consider $u_i/u_j-v_i/v_j$, and apply mean value theorem
to function $g$, defined as $g(x, y)=x/y$, we show that there exists
a positive constant $D_2$ such that
$$
d_2(u, v) < D_2 d_{\mathbf{E}}(u, v).
$$
Namely $d_2 \sim d_E$. Thus the claim in this Proposition follows,
namely there exist two positive constant $C_1 < C_2$ such that for
any two points $u, v \in S$,
$$
C_1 d_\textbf{B}(u, v) < d_\textbf{E}(u, v) < C_2 d_\textbf{B}(u,
v).
$$
\end{proof}

\section{Proof of Lemma~\ref{bounds}:}
\label{lem-bounds}

Recall that for a non-negative matrix $B$, the \emph{canonical form}
of $B$ is:
$$
B=\left [  \begin{array}{cccc}
            B_{11}& B_{12} & \cdots & B_{1n}\\
            0& B_{22} & \cdots & B_{2n}\\
            \vdots & \vdots & \ddots & \vdots\\
            0& 0& \cdots & B_{nn}\\
  \end{array} \right ],
$$
where $B_{ii}$ is either an irreducible matrix (called {\em
irreducible components}) or a $1 \times 1$ zero matrix.

Condition 2 in Theorem~\ref{unambig-bdry} is equivalent to the
statement that $B = B(\vec\eps_0)$ has a unique irreducible
component of maximal spectral radius and that this component is
primitive. Let $C$ denote the square matrix obtained by restricting
$B$ to this component and let $S_C$ denote the set of indices
corresponding to this component. Let $\lambda_1$ denote the spectral
radius of $B$, equivalently the spectral radius of $C$.

Let $\lambda_1(\vec\eps)$ denote the largest, in modulus, eigenvalue
of $B(\vec\eps)$.  Since the entries of $B(\vec\eps)$ are analytic
in $\vec\eps$ and $\lambda_1$ is simple, it follows that if the
complex neighborhood $\Omega$ is chosen sufficiently small, then
$\lambda_1(\vec\eps)$ is analytic function of $\vec\eps\in \Omega$.

The columns (resp., rows) of $Adj(\lambda_1(\vec\eps)I -
B(\vec\eps))$ are right (resp., left) eigenvectors of $B(\vec\eps)$
corresponding to $\lambda_1(\vec\eps)$.  By choosing $x(\vec\eps)$
(resp. $y(\vec\eps)$) to be a fixed column (resp. row) of
$Adj(\lambda_1(\vec\eps)I - B(\vec\eps))$ and then replacing
$x(\vec\eps)$ and $y(\vec\eps)$ by appropriately rescaled versions,
we may assume that:

\begin{itemize}
\item $x(\vec\eps_0), y(\vec\eps_0) \geq 0$, and they are positive on $S_C$
\item $y(\vec\eps) \cdot x(\vec\eps) =1$
\item $x(\vec\eps)$ and $y(\vec\eps)$ are analytic in $\vec\eps \in \Omega$
\end{itemize}

Let $$V(\vec\eps) = \lambda_1(\vec\eps)  x(\vec\eps) \cdot
y(\vec\eps)$$ and
$$U(\vec\eps) = B(\vec\eps) - V(\vec\eps).$$
Then $V(\vec\eps)$ is the restriction of $B(\vec\eps)$ to the
subspace corresponding to $\lambda_1(\vec\eps)$ and $U(\vec\eps)$ is
the restriction to the subspace corresponding to the remainder of
the spectrum of $B(\vec\eps)$. It follows that
$$U(\vec\eps)V(\vec\eps) = 0 = V(\vec\eps) U(\vec\eps).$$

Let $\mu(\vec\eps)$ denote the spectral radius of $U(\vec\eps)$. By
condition $2$, $\mu(\vec\eps_0) < \lambda_1(\vec\eps_0)$.  Thus,
there is a constant $\nu
> 0$ such that if the neigbourhood $\Omega$ is sufficiently small,
then for all $\vec\eps \in \Omega$
$$
\mu(\vec\eps) < \nu <  |\lambda_1(\vec\eps)|.
$$
Thus, by Lemma~\ref{Senata}, and making still $\Omega$  smaller if
necessary, there is a constant $K_1>0 $ such that for all $i,j$, all
$n$ and all $\vec\eps \in \Omega$,
\begin{equation}
\label{ULK} |U_{ij}^n(\vec\eps)| < K_1 \nu^n.
\end{equation}

Let $r = r(\vec\eps_0)$, $c = c(\vec\eps_0)$, $x=x(\vec\eps_0)$ and
$y=y(\vec\eps_0)$. In the following we will show that the
irreducibility of $\Delta$ will rule out the possibility that $c$ is
non-zero only in non-maximal spectral radius irreducible components
of $B$, and so we can extend $a(\vec{\eps}, n)$ and $b(\vec{\eps},
n)$ from real to complex.

Let $s_0 \in S_C$. Since $\Delta(\vec\eps_0)$ is irreducible and $r$
is nonnegative, but not the zero vector, for some $j_0$,
$(rB^{j_0})_{s_0} > 0$. Similarly, for any index $s_1$ other than 1
of the underlying Markov chain, there exists $j_1$ such that
$B^{j_1}_{s_0s_1}
> 0$. Choose $s_1$ to be any index such that $c_{s_1} > 0$. Since $C$ is primitive, it
then follows that there is a constant $K_2$ such that for
sufficiently large $n$,
$$
r x \cdot y c \lambda_1^n + r U^n c=r V^n c + r U^n c = r B^n c  >
K_2 \lambda_1^n,
$$
which by (\ref{ULK}) implies that $r x \cdot y c > 0$. Therefore if
$\Omega$ is sufficiently small, there exists a positive constant
$K_4$ such that
$$
|r(\vec\eps) x(\vec\eps) \cdot y(\vec\eps) c(\vec\eps)| > K_4,
$$
for $\vec\eps \in \Omega$.

Let $K_3$ be an upper bound on the entries of $|x(\vec\eps)|,
|y(\vec\eps)|, |r(\vec\eps)|$ and $|c(\vec\eps)|$.

Thus, for all $n$ and all $\vec\eps \in \Omega$, we have
$$
|r(\vec\eps) B^n (\vec\eps) c(\vec\eps)| \le | r(\vec\eps) U^n
(\vec\eps) c(\vec\eps)| +|r(\vec\eps) V^n (\vec\eps)  c(\vec\eps)|
 \le |{\cal B}|^2 K_3^2 K_1 \nu^n +|{\cal B}|^2 K_3^4
|\lambda_1(\vec\eps)|^{n }
$$
and
$$
|r(\vec\eps) B^n (\vec\eps) c(\vec\eps)| \ge |r(\vec\eps) V^n
(\vec\eps) c(\vec\eps)| - | r(\vec\eps) U^n (\vec\eps) c(\vec\eps)|
\ge K_4 |\lambda_1(\vec\eps)|^n - |{\cal B}|^2 K_3^2 K_1 \nu^n.
$$
With similar upper and lower bounds for $|(r(\vec\eps) B^n
(\vec\eps) {\bf 1}|$, it follows that for sufficiently large $n$ and
all $\vec\eps \in\Omega$,
$$
\frac{\pi_1(\vec\eps)r(\vec\eps)B(\vec\eps)^n {\bf
1}}{\pi_1(\vec\eps)r(\vec\eps)B(\vec\eps)^{n-1} {\bf 1}}
$$
and
$$
\frac{\pi_1(\vec\eps)r(\vec\eps)B(\vec\eps)^{n-1}
c(\vec\eps)}{\pi_1(\vec\eps)r(\vec\eps)B(\vec\eps)^{n-1} {\bf 1}}
$$
are uniformly bounded from above and away from zero. By condition
$1$, for any finite collection of $n$, there is a (possibly smaller)
neighborhood $\Omega$ of $\vec{\varepsilon}_0$, such that for all
$\vec\eps \in \Omega$, these quantities are uniformly bounded from
above and away from zero. This completes the proof of
Lemma~\ref{bounds} ( and therefore the proof of sufficiency for
Theorem~\ref{unambig-bdry}.)

\section{$\vec{\varepsilon} \mapsto \nu^{\vec{\varepsilon}}$ is analytic}  \label{equilibrium}
In this appendix, we follow the notation in Section~\ref{xxx-II}.
Let $\tau: \mathcal{X} \rightarrow \mathcal{X}$ be the right shift
operator, which is a continuous mapping on $\mathcal{X}$ under the
topology induced by the metric $d$. For $f \in C(\mathcal{X})$, one
defines the {\em pressure}  via a variational principle~\cite{ru78}:
$$
P(f)=\sup_{\mu \in M(\mathcal{X}, \tau)} \left( H_{\mu}(\tau)+\int f
d\mu \right),
$$
where $M(\mathcal{X}, \tau)$ denotes the set of $\tau$-invariant
probability measures on $\mathcal{X}$ and $H_{\mu}(\tau)$ denotes
measure-theoretic entropy. A member $\mu$ of $M(\mathcal{X}, \tau)$
is called an \emph{equilibrium state} for $f$ if $P(f)=H_{\mu}
(T)+\int f d\mu$.

For $f \in C(\mathcal{X})$ the Ruelle operator $\mathcal{L}_f :
C(\mathcal{X}) \rightarrow C(\mathcal{X})$ is defined~\cite{ru78} by
$$
(\mathcal{L}_f h)(x)=\sum_{y \in \tau^{-1}x} e^{f(y)} h(y).
$$

The connection between pressure and the Ruelle operator is as
follows~\cite{ru78, ru97}. When $f \in F^{\theta}$, $P(f)$  is $\log
\lambda$, where $\lambda$ is the spectral radius of $\mathcal{L}_f$.
The restriction of $\mathcal{L}_f$ to $F^{\theta}$ still has
spectral radius $\lambda$, and $\lambda$ is isolated from all other
eigenvalues of the restricted operator. Using this, Ruelle applied
standard perturbation theory for linear operators~\cite{ka76} to
conclude that pressure $P(f)$ is real analytic on $F^{\theta}$.
Moreover, he showed that each $f \in F^{\theta}$ has a unique
equilibrium state $\mu_f$ and the first order derivative of $f
\mapsto P(f)$ on $F^{\theta}$ is $\mu_f$, viewed as a linear
functional on $F^{\theta}$. So, the analyticity of $P(f)$ implies
that the equilibrium state $\mu_f$ is also analytic in $f \in
F^{\theta}$.

We first claim that for $f(\vec{\varepsilon},z) = \log
p^{\vec{\varepsilon}}(z_0|z_{-\infty}^{-1})$, we have
$\mu_{f(\vec{\varepsilon},\cdot)} = \nu^{\vec{\varepsilon}}$ as in
(\ref{ppp}).

To see this, first observe that the spectral radius $\lambda$ of
$\mathcal{L} = \mathcal{L}_{f(\vec{\varepsilon},\cdot)}$ is 1; this
follows from the observations:
\begin{itemize}
\item the function $\bar 1$ which is identically
1 on $\mathcal{X}$ is a fixed point of $\mathcal{L}$ ~~~ -- and --
\item (see Proposition 5.16 of~\cite{ru78}) $\mathcal{L}^n(\bar
1)/\lambda^n$ converges to a strictly positive function.
\end{itemize}
Thus $P({f(\vec{\varepsilon},\cdot)}) = 0$.  So, for
$\mu^{\vec{\varepsilon}} = \mu_{f(\vec{\varepsilon}, \cdot)}$, we
have
$$ h_{\mu^{\vec{\varepsilon}}}(\tau)+\int f(\vec{\varepsilon},\cdot)
d\mu^{\vec{\varepsilon}} = 0.$$ But from (\ref{ent-form}), we have
$$ h_\nu(\tau)+\int f(\vec{\varepsilon},\cdot)
d\nu^{\vec{\varepsilon}} = 0.$$ By uniqueness of the equilibrium
state, we thus obtain $\mu_{f(\vec{\varepsilon}, \cdot)} =
\nu^{\vec{\varepsilon}}$  as claimed.

Since $\vec{\varepsilon} \mapsto f(\vec{\varepsilon},\cdot)$ is
analytic, it then follows that  $\vec{\varepsilon} \mapsto
\nu^{\vec{\varepsilon}}$ is analytic, thereby completing the proof
of Theorem~\ref{main-1}.


\begin{thebibliography}{10}

\bibitem{ar94a}
L.~Arnold, V.~M.~Gundlach and L.~Demetrius.
\newblock Evolutionary formalism for products of positive random matrices.
\newblock {\em Annals of Applied Probability}, 4:859--901, 1994.

\bibitem{bi62}
J.~J.~Birch.
\newblock Approximations for the entropy for functions of {M}arkov chains.
\newblock {\em Ann. Math. Statist.}, 33:930--938, 1962.

\bibitem{bl57}
D.~Blackwell.
\newblock The entropy of functions of finite-state {M}arkov chains.
\newblock {\em Trans. First Prague Conf. Information Thoery, Statistical
  Decision Functions, Random Processes}, pages 13--20, 1957.

\bibitem{ca81}
M.~Cassandro and E.~Olivieri.
\newblock Renormalization group and
analyticity in one dimension: A proof of Dobrushin's theorem
\newblock {\em Commun. Math. Phys.}, 80, 255-269, 1981.

\bibitem{ch03}
J.~R.~Chazottes and E.~Ugalde.
\newblock Projection of {M}arkov measures may be Gibbsian.
\newblock {\em J. Statist. Phys.}, Volume 111, Numbers 5-6, 1245-1272.

\bibitem{do73}
R.~L.~Dobrushin.
\newblock Analyticity of correlation functions in one-dimensional classical
systems with slowly decreasing potentials.
\newblock Commun. Math. Phys. 32, 269-289, 1973.

\bibitem{eg04}
S.~Egner, V.~Balakirsky, L.~Tolhuizen, S.~Baggen and H.~Hollmann.
\newblock On the entropy rate of a hidden {M}arkov model.
\newblock In {\em Proceedings of the 2004 IEEE International Symposium on
  Information Theory}, page~12, Chicago, U.S.A., 2004.

\bibitem{gm05}
G.~Han and B.~Marcus.
\newblock Analyticity of entropy rate of a hidden {M}arkov chain
\newblock In Proc. of IEEE International Symposium on Information Theory,
Adelaide, Australia, September 4-September 9 2005, pages 2193-2197.


\bibitem{gh95u}
R.~Gharavi and V.~Anantharam.
\newblock An upper bound for the largest {L}yapunov exponent of a {M}arkovian
  product of nonnegative matrices.
\newblock Preprint, Janurary 1995.

\bibitem{ho03}
T.~Holliday, A.~Goldsmith and P.~Glynn.
\newblock On entropy and {L}yapunov exponents for finite state
channels. 2003. Available at
http://wsl.stanford.edu/Publications/THolliday/Lyapunov.pdf.

\bibitem{ja04}
P.~Jacquet, G.~seroussi and W.~Szpankowski.
\newblock On the entropy of a hidden {M}arkov process.
\newblock In {\em Proceedings of the 2004 IEEE International Symposium on
  Information Theory}, page~10, Chicago, U.S.A., 2004.

\bibitem{ka76}
T.~Kato.
\newblock {\em Perturbation Theory for Linear Operators}.
\newblock Springer Verlag, Berlin-Heidelberg-New York, 1976.

\bibitem{lm95}
D.~Lind and B.~Marcus.
\newblock {\em An Introduction to Symbolic Dynamics and Coding}.
\newblock Cambridge University Press, 1995.

\bibitem{lo98}
J. L\"{o}rinczi, C. Maes and K. V. Velde.
\newblock Transformations of Gibbs measures.
\newblock {\em Probab. Theory Relat. Fields}, Volume 112, 121-147,
1998.

\bibitem{ma84a1}
B.~Marcus, K.~Petersen and S.~Williams.
\newblock Transmission rates and factors of {M}arkov chains.
\newblock {\em Contemporary Mathematics}, 26:279--294, 1984.

\bibitem{mu78}
A.~Mukherjea and K.~Pothoven.
\newblock {\em Real and functional analysis}.
\newblock Plenum Press, New York, 1978.

\bibitem{na81}
L.~Nachbin.
\newblock {\em Introduction to functional analysis : Banach
spaces and differential calculus}.
\newblock New York : M. Dekker, 1981.

\bibitem{on93}
A.~Onishchik.
\newblock {\em Lie groups and Lie algebra I}.
\newblock Encyclopaedia of mathematical sciences ; v. 20. Springer-Verlag, 1993.

\bibitem{or03}
E.~Ordentlich and T.~Weissman.
\newblock On the optimality of symbol by symbol filtering and denoising.
\newblock {\em Information Theory, IEEE Transactions}, Volume 52,  Issue 1,  Jan. 2006 Page(s):19 -
40.

\bibitem{or04}
E.~Ordentlich and T.~Weissman.
\newblock New bounds on the entropy rate of hidden {M}arkov process.
\newblock {\em Information Theory Workshop}, 2004. IEEE 24-29 Oct. 2004 Page(s):117
- 122

\bibitem{pe90}
Y.~Peres.
\newblock {\em Analytic dependence of Lyapunov exponents on transition
probabilities}, volume 1486 of {\em Lecture Notes in Mathematics,
Lyapunov's exponents, {P}roceedings of a {W}orkshop}.
\newblock Springer Verlag, 1990.

\bibitem{pe92}
Y.~Peres.
\newblock Domains of analytic continuation for the top {L}yapunov exponent.
\newblock {\em Ann. Inst. H. Poincar\'e Probab. Statist.}, 28(1):131--148,
  1992.

\bibitem{pe03}
K.~Petersen, A.~Quas and S.~Shin.
\newblock Measures of maximal relative entropy.
\newblock {\em Ergod. Th. and Dynam. Sys.}, 23, 207-223, 2003

\bibitem{ru78}
D.~Ruelle.
\newblock {\em Thermodynamic formalism : the mathematical structures of
  classical equilibrium statistical mechanics}.
\newblock Addison-Wesley Pub. Co., Advanced Book Program, Reading, Mass, 1978.

\bibitem{ru79}
D.~Ruelle.
\newblock Analyticity properties of the characteristic exponents of random
  matrix products.
\newblock {\em Adv. Math.}, 32:68--80, 1979.

\bibitem{ru97}
D.~Ruelle.
\newblock Differentiation of {SRB} states.
\newblock {\em Comm. Math. Phys.}, 187(1):227--241, 1997.

\bibitem{se80}
E.~Seneta.
\newblock {\em Springer Series in Statistics}.
\newblock Non-negative Matrices and Markov Chains. Springer-Verlag, New York
  Heidelberg Berlin, 1980.

\bibitem{sh92}
B.~V.~Shabat.
\newblock {\em Introduction to complex analysis}.
\newblock Translations of mathematical monographs ; v. 110. American
  Mathematical Society, Providence, R.I., 1992.


\bibitem{ta02}
J.~L.~Taylor.
\newblock {\em Several complex variables with connections to algebraic geometry
and Lie groups}.
\newblock American Mathematical Society, Providence, R.I., 2002.

\bibitem{wa82}
P.~Walters.
\newblock {\em An introduction to ergodic theory}. volume~79 of {\em Graduate
  texts in mathematics}.
\newblock Springer-Verlag, New York, 1982.

\bibitem{yo74}
K.~Yosida.
\newblock {\em Functional analysis}, 4th edition.
\newblock Springer-Verlag, Berlin, 1974.

\bibitem{zu04}
O.~Zuk, I.~Kanter and E.~Domany.
\newblock Asymptotics of the entropy rate for a hidden {M}arkov process.
\newblock {\em J. Stat. Phys.}, 121(3-4): 343-360 (2005)

\bibitem{zu05}
O.~Zuk, E.~Domany, I.~Kanter, and M.~Aizenman.
\newblock Taylor series expansions for the entropy rate of
Hidden Markov Processes.
\newblock ICC 2006, Istanbul.

\end{thebibliography}

\end{document}